\documentclass[12pt, twoside, a4paper, fleqn]{article}

\usepackage{latexsym}
\usepackage{amscd}
\usepackage{theorem}
\usepackage{pifont}
\usepackage{mathbbol}
\usepackage{amsfonts}
\usepackage{xspace}
\usepackage{amssymb}
\usepackage{fancyhdr}
\usepackage{mathrsfs}
\usepackage{amsmath}%
\usepackage{graphics}%
\usepackage{graphicx}%
\usepackage{euscript}%
\usepackage{psfrag}%
\usepackage{empheq}%
\usepackage{upgreek}
\usepackage{eufrak}

\DeclareMathAlphabet{\mathpzc}{OT1}{pzc}{m}{it}

{\theorembodyfont{\slshape} \newtheorem{theorem}{Theorem}[section]}
{\theorembodyfont{\slshape} \newtheorem{lemma}[theorem]{Lemma}}
{\theorembodyfont{\slshape} \newtheorem{proposition}[theorem]{Proposition}}
{\theorembodyfont{\slshape} }
{\theorembodyfont{\slshape} \newtheorem{remark}[theorem]{Remark}}
{\theorembodyfont{\slshape} \newtheorem{claim}[theorem]{Claim}}

\newenvironment{proof}[1][\proofname]{\noindent \textbf{Proof of {#1}:\ }}{$\hfill{\Box}$}
\newenvironment{prooof}{\noindent\textbf{Proof:\ }}{$\hfill{\Box}$}

\numberwithin{equation}{section}

\title{\textsc{Simultaneous Action of Finitely Many Interval Maps: Some Dynamical and Statistical Properties}}                

\author{Aswin Gopakumar \\ {\tt aswin15@iisertvm.ac.in} \bigskip \\ Kirthana Rajasekar \\ {\tt kirthanarajasekar15@iisertvm.ac.in} \bigskip \\ Shrihari Sridharan \\ {\tt shrihari@iisertvm.ac.in} \bigskip 
 \bigskip \\ {\sl Indian Institute of Science Education and Research} \\ {\sl Thiruvananthapuram (IISER-TVM), India.}} 

\date{\today}

\begin{document}

\maketitle
\thispagestyle{empty}

\begin{abstract} 
\noindent 
In this paper, we consider finitely many interval maps simultaneously acting on the unit interval $I = [0,\, 1]$ in the real line $\mathbb{R}$; each with utmost finitely many jump discontinuities and study certain important statistical properties. Even though we use the symbolic space on $N$ letters to reduce the case of simultaneous dynamics to maps on an appropriate space, our aim in this paper remains to resolve ergodicity, rates of recurrence, decay of correlations and invariance principles leading upto the central limit theorem for the dynamics that evolves through simultaneous action. In order to achieve our ends, we define various Ruelle operators, normalise them by various means and exploit their spectra. 
\end{abstract}
\bigskip \bigskip

\begin{tabular}{l c l}
\textbf{Keywords} & : & Growth of typical trajectories; \\ 
& & Invariance principles; \\ 
& & Ruelle operator and the pressure function; \\ 
& & Simultaneous action of finitely many interval maps. \\
& & \\
\textbf{AMS Subject} & : & 37E05, 37C35, 37D35, 37B10. \\
\textbf{Classifications} & & 
\end{tabular}
\bigskip 
\bigskip 

\newpage 

\section{Introduction}

\noindent 
Various dynamical properties and statistical properties help us understand the behaviour of dynamics caused by the action of a transformation $T$ on some phase space $X$. Important among such properties include the Birkhoff's pointwise ergodic theorem, asymptotic estimates on rates of recurrence of typical orbits, decay of correlations, invariance principles, central limit theorem, law of iterated logarithms {\it etc}. Each of these theorems provide us a deeper glimpse into the structure, the dynamics builds in its phase space of action or an invariant subset, thereof. 
\medskip 

\noindent 
Birkhoff's pointwise ergodic theorem observes a considered dynamical system through a real-valued continuous function and states that the sequence of local time averages along the orbit of any typical point converges to the global space average, whenever the transformation $T$ acting on $X$ is ergodic. Though the result is extremely strong, it does allow some points (though negligible, meaning with collective measure zero) to fluctuate from this mean behaviour. Thus, an interesting study in the dynamics of such ergodic systems is to obtain a good understanding of the set of points that violate the Birkhoff's ergodic conditions. An easy way to approach this subject locally is to work out the ergodic sums of the observables and consider the cardinality of the set of points whose ergodic sum calculated at various times remain inside some chosen interval $[a,\, b] \subset \mathbb{R}$. However, on a global scale, an alternate way to understand the deviation from the average behaviour of typical orbits is achieved by formalising the central limit theorem. 
\medskip 

\noindent 
The central limit theorem is an important tool in mathematics that distributes the random variables along a bell-curve (normal distribution), as more and more independent random variables are appropriately included under the ambit of study. This is a central object of investigation in understanding deterministic dynamical systems, owing to its natural appeal, when we consider the various orbits in the phase space. However, as important as the central limit theorem is, we see that they are subsumed by more general invariance principles. A sequence of random variables $\big\{ X_{n} \big\}_{n\, \ge\, 1}$ is said to satisfy an almost sure invariance principle if the sequence can be approximated almost surely by another sequence, preferably with certain desired properties and with a relatively small margin of error. 
\medskip 

\noindent 
Several mathematicians have studied these properties in many deterministic dynamical systems, where the phase space is a compact interval of the real line, \cite{md:86, ci:96, lsv:99, ps:02, cr:07, cm:15}, the Julia set of some rational map that occurs as a compact subset of the Riemann sphere, \cite{du:91, dpu:96, ss:07, ss:09}, {\it etc.}, however, with a single transformation acting on the appropriate space, based on which the dynamics evolves. Examples of continuous time dynamical systems that has remained in the focus of the dynamics community include expanding flows restricted on a compact subset of the Riemannian manifold, \cite{dp:84, spl:89, mp:91} or some mixing Axiom $A$ diffeomorphism restricted on a basic set, \cite{mr:73, ks:90, rs:93, na:00, ps:01}. A particularly desirable feature of all the above-mentioned maps restricted on their respective sets is that they can be studied through an associated symbolic model \cite{rb:73, mr:73}. There are also various studies carried out by several mathematicians that analyse statistical results in various settings of dynamical systems. Prominent among them include \cite{ps:75, hh:80, cp:90, pp:90, ps:94, lsy:99, si:99, fmt:03, mtk:05, mn:05, hntv:17}. 
\medskip 

\noindent 
What we intend to investigate in this paper is slightly richer in dynamics, than what is explained so far. In this paper, we consider the compact unit interval $[0,\, 1] \subset \mathbb{R}$; however with finitely many maps acting on the space simultaneously. Thus, the dynamics evolves along the multiple branches provided by each of these maps. In fact, we work with finitely many interval maps defined on $[0,\, 1]$; each of which has a discrete set of utmost finitely many discontinuities. As an expert reader may realise, these results are readily transferable to various settings including the simultaneous action of finitely many rational maps restricted on the appropriate Julia set, as defined in \cite{hs:00} or to the action of a holomorphic correspondence restricted on the support of its Dinh-Sibony measure, as defined in \cite{bs:16}. We shall explain our claim of transferability of the main theorems of this paper, as written above, in the final section, \eqref{concl}. 
\medskip 

\noindent 
This paper is structured as follows: In the next two introductory sections \eqref{prelims} and \eqref{peri}, we narrate the basic settings of this paper and develop certain notations and dynamical notions, however only as skeletal to enable us to state the main theorems of this paper in section \eqref{main}. The main results of this paper describe the ergodicity of the system in theorems \eqref{erg1} and \eqref{erg2}, the rates of recurrence in theorems \eqref{ror1} and \eqref{ror2}, the exponential decay of correlations in theorems \eqref{doc1} and \eqref{doc2}, almost sure invariance principles in theorems \eqref{asip1} and \eqref{asip2} and a few more statistical properties such as the central limit theorems and the laws of iterated logarithms in theorems \eqref{osp1} and \eqref{osp2}. The reason why each theorem appears twice in the list above will be clear, by the time we reach section \eqref{main}. In section \eqref{SigmaNplus}, we recall the setting of symbolic dynamics that comes in handy as a book-keeping mechanism in our study. In sections \eqref{vro}, \eqref{spectrum} and \eqref{Normalising}, we define three kinds of Ruelle operators on the appropriate Banach space of continuous functions and H\"{o}lder continuous functions defined on the phase spaces that interest us, compare their spectra and normalise them in different ways in order that they help us in proving our main theorems. Having achieved these, we embark on writing the proofs of the main theorems in sections \eqref{ergproof}, \eqref{rorproof}, \eqref{docsec}, \eqref{asipsec} and \eqref{seccor}. We conclude the paper with a few remarks in section \eqref{concl}. 
\bigskip 

\section{Preliminaries and the pressure function} 
\label{prelims} 

\noindent 
In this section, we explain the setting of our paper and define certain basic terminologies that help us in constructing the necessary notions to state our main results. 
\medskip 

\noindent 
Let $I$ denote the unit interval on the real line, {\it i.e.}, $I = [0,\, 1]$. We are interested in studying the dynamics of finitely many interval maps acting simultaneously on $I$, \textit{i.e.}, given $N \in \mathbb{N}$ and $1 \leq d \leq N$, we consider the interval maps $T_{d} : I \longrightarrow I$ of degree $(d + 1)$ given by 
\[ T_{d}\, (x)\ \ :=\ \ (d + 1)\; x \pmod 1. \] 
The simultaneous action is explained as follows: For any $x_{0} \in I$, its forward orbit at times $t = 0,\, 1,\, 2,\, \cdots,\, n,\, \cdots$ is defined as 
\begin{equation} 
\label{1storbit}
\left\{ x_{0},\, x_{1} \in \bigcup_{d\, =\, 1}^{N} T_{d} (x_{0}),\, x_{2} \in \bigcup_{d\, =\, 1}^{N} T_{d} (x_{1}),\, \cdots,\, x_{n} \in \bigcup_{d\, =\, 1}^{N} T_{d} (x_{n - 1}), \cdots \right\}.  
\end{equation}
Thus, at every stage, we have $N$ many maps to choose from to move forward and the totality of all these branches describe the forward orbit. Observe that the dynamics that arises out of such a process can also be described by the action of a semigroup generated by the same interval maps, $\mathscr{S} = \big\langle T_{1},\, T_{2},\, \cdots,\, T_{N} \big\rangle$, or as a correspondence on $I \times I$, appropriately defined. 
\medskip 

\noindent 
Suppose $\mathcal{C} (V, \mathbb{F})$ denotes the space of all continuous functions defined on the space $V$ that takes values on the field $\mathbb{F}$. Then, for any $f \in \mathcal{C}(I, \mathbb{C})$, the set of all complex-valued continuous functions defined on $I$, we define a \emph{composition operator}, $\mathscr{O}$ by 
\[ \mathscr{O} (f) \in \bigcup\limits_{d\, =\, 1}^{N} \big( f \circ T_{d} \big). \] 
Along every orbit of the point $x_{0}$, as described in \eqref{1storbit}, this composition operator chooses the map $T_{d}$ every time in such a fashion that $T_{d} (x_{k - 1})\, =\, x_{k}$. Hence, even though $\mathscr{O} (f)$ is not single-valued, we have by definition that $\left( \mathscr{O} (f) \right) (x_{0})$ to be single-valued along every chosen orbit of $x_{0}$. We further describe this idea in section \eqref{peri}. 
\medskip 

\noindent 
To assist us in this study, we will make use of the space consisting of infinitely long words on $N$ symbols, {\it i.e.}, suppose $S = \big\{ 1,\, 2,\, \cdots,\, N \big\}$, we consider 
\[ \Sigma_{N}^{+}\ \ :=\ \ S^{\mathbb{Z}_{+}}\ \ =\ \ \Big\{ w = \left( w_{1}\, w_{2}\, \cdots\, w_{n}\, \cdots \right) : w_{i} \in S \Big\}. \] 
As we shall explain in section \eqref{SigmaNplus}, $\Sigma_{N}^{+}$ is a compact measurable metric space equipped with the Bernoulli measure, where the shift map $\sigma$ defined by $\left( \sigma (w) \right)_{n} = w_{n\, +\, 1}$ is continuous and non-invertible, but a local homeomorphism. 
\medskip 

\noindent 
From now on, we denote by $X$ the product phase space given by $X := \Sigma_{N}^{+} \times I$, where we define a skew-product map $T$ as 
\begin{equation} 
\label{spm} 
T (w,\, x)\ \ :=\ \ ( \sigma w,\, T_{w_{1}} x),\ \ \ \ \text{where}\ \ w = (w_{1}\, w_{2}\, \cdots). 
\end{equation} 

\noindent 
By a standard argument due to Tychonoff, as may be found in \cite{mj:00}, we consider the natural product topology on $X$ that gives rise to the metric $d_{X} (\cdot,\, \cdot)$ on $X$. Further, we also have the product sigma-algebra and the product measure defined on $X$. Let $\mu$ denote some $T$-invariant measure supported on $X$. For example, the appropriate product measure of the Bernoulli measure on cylinder sets of $\Sigma_{N}^{+}$ and the Lebesgue measure on open intervals of $I$ is a $T$-invariant probability measure on $X$. 
\medskip 

\noindent 
The definition of the skew-product map $T$ entails that the forward orbit of $(w, x)$ at times $t = 0,\, 1,\, 2,\, \cdots,\, n,\, \cdots$ under $T$ is given by 
\[ \Big\{ (w,\, x),\; (\sigma w,\, T_{w_{1}} x),\; \left( \sigma^{2} w,\, \left( T_{w_{2}} \circ T_{w_{1}} \right) x \right),\; \cdots,\; \left( \sigma^{n} w,\, \left( T_{w_{n}} \circ T_{w_{n - 1}} \cdots \circ T_{w_{1}} \right) x \right),\; \cdots \Big\}.  \] 
Thus, for a chosen $w$ in $\Sigma_{N}^{+}$, the natural projection on the second co-ordinate $\mathbf{Proj}_{2} : X \longrightarrow I$ captures the sectional idea behind the orbit of $x \in I$ as described in \eqref{1storbit}. Taking the union over all possible $w \in \Sigma_{N}^{+}$ captures the idea in its entirety. 
\medskip 

\noindent 
For ease of notations, we fix little letters like $f,\, g,\, h,$ {\it etc}. to represent functions defined on the interval, $I$ and use big letters like $F,\, G,\, H$ {\it etc}. to represent functions defined on the product space $X = \Sigma_{N}^{+} \times I$. Although one may think of $f$ as being a restriction of $F$ on the interval, {\it i.e.}, $F = f \circ \mathbf{Proj}_{2}$ that yields $F((w,\, x)) = f(x)$, it need not be the case always. 
\medskip 

\noindent 
The space $X$ now facilitates us to redefine the composition operator, in this setting denoted by $\mathscr{Q}$ defined on $\mathcal{C} (X, \mathbb{C})$ given by  
\[ \mathscr{Q} (F)\ \ :=\ \ F \circ T. \] 

\noindent 
Let $\mathscr{F}_{\alpha}(X, \mathbb{C})$ denote the space of all complex-valued $\alpha$-H\"{o}lder continuous functions defined on $X$. For $F \in \mathscr{F}_{\alpha}(X, \mathbb{C})$, we define the following norm, 
\[ \big\Vert F \big\Vert_{\alpha}\ \ :=\ \ \big\vert F \big\vert_{\alpha} + \big\Vert F \big\Vert_{\infty}, \] 
where 
\[ \big\vert F \big\vert_{\alpha}\ \ :=\ \ \sup\limits_{(w,\, x)\; \neq\; (v,\, y)} \left\{ \frac{\big\vert F((w,\, x))\; -\; F((v,\, y)) \big\vert}{\big( d_{X}((w,\, x),\; (v,\, y)) \big)^{\alpha}}\ :\ (w,\, x),\; (v,\, y) \in X \right\} \] 
denotes the $\alpha$-H\"older semi-norm and  $\big\Vert F \big\Vert_{\infty}$ denotes the usual supremum norm. Then, $\mathscr{F}_{\alpha}(X, \mathbb{C})$ is a Banach space under the norm $\big\Vert \cdot \big\Vert_{\alpha}$. 
\medskip 

\noindent 
Given any function $F \in \mathcal{C}(X, \mathbb{R})$, its \emph{pressure} is defined as 
\begin{equation} 
\label{pressure} 
\mathfrak{P} (F)\ \ :=\ \ \sup \left\{ h_{\mu} (T) + \int\! F d \mu \right\}, 
\end{equation} 
where the supremum is taken over all $T$-invariant probability measures supported on $X$. Further, $h_{\mu} (T)$ is the measure theoretic entropy of $T$ with respect to $\mu$. An \emph{equilibrium measure} for the function $F$ denoted by $\mu_{F}$ is defined as that measure for which the supremum is attained in the definition of pressure, as stated in \eqref{pressure}. The unique existence of $\mu_{F}$ for every $F \in \mathscr{F}_{\alpha} (X, \mathbb{R})$ is assured by Denker and Urbanski in \cite{du:91} and Sumi and Urbanski in \cite{su:09}, for an analogous setting. 
\medskip

\noindent 
The pressure function and the equilibrium measure have respective analogues for dynamics under simultaneous actions of the interval maps. We will establish the same later in section \eqref{Normalising}. However, for the sake of stating the results, we mention the following properties. Given $f \in \mathscr{F}_{\alpha}(I, \mathbb{R})$, its \textit{pressure} under simultaneous dynamics that we denote by $\mathbb{P}(f)$ coincides with the quantity $\mathfrak{P}(f \circ \mathbf{Proj}_{2})$. Similarly, the measure $\mathfrak{m}_{f}$ on $I$, whose relation with $\mu_{f \circ \mathbf{Proj}_{2}}$, which will be an easy observation once defined, is given by,
\[\int\! g\, d \mathfrak{m}_{f}\ \ =\ \ \int\! ( g \circ \mathbf{Proj}_{2})\, d \mu_{f\, \circ\, \mathbf{Proj}_{2}},\ \ \ \forall g \in \mathscr{F}_{\alpha}(I, \mathbb{R}).\] 
\bigskip 

\section{Periodic points and other dynamical notions}
\label{peri}

\noindent 
A point $(w,\, x) \in X$ is said to be a \emph{periodic point} of period $p$ for $T$ iff $\sigma^{p} w\; =\; w$ and $\left( T_{w_{p}} \circ T_{w_{p - 1}} \circ \cdots \circ T_{w_{1}} \right)\, x\; =\; x$. We denote the set of all $p$-periodic points by ${\rm Fix}_{p} (T)$. Once we determine the periodic points of $T$, it becomes easier for us to identify the periodic orbits of simultaneous action of the $N$ many interval maps. 
\medskip 

\noindent 
We first introduce a few notations here. For any $x_{1} \in I$, let $\mathscr{R}_{n} (x_{1})$ denote the set of all rays starting from $x_{1}$ that describe the initial $n$-long itinerary of the trajectory of $x_{1}$ in the order that the point visits, {\it i.e.}, 
\[ \mathscr{R}_{n} (x_{1})\ \ :=\ \ \Big\{ \left( x_{1},\, x_{2},\, \cdots,\, x_{n} \right) \in I^{n}\ :\ \forall 2 \le k \le n,\ \exists 1 \le d \le N\ \text{such that}\ T_{d} (x_{k - 1}) = x_{k} \Big\}. \] 
Further, we denote by $\mathscr{R} (x_{1}) = \mathscr{R}_{\infty} (x_{1})$ the set of all infinite rays starting from $x_{1}$, where we produce each point in $\mathscr{R}_{n} (x_{1})$ to an infinitely long sequence, as allowed by the dynamics, {\it i.e.},  
\[ \mathscr{R} (x_{1})\ \ :=\ \ \Big\{ \left( x_{1},\, x_{2},\, \cdots\, \right) \in I^{\mathbb{Z}_{+}}\ :\ \forall k \ge 2,\ \exists 1 \le d \le N\ \text{for which}\ T_{d} (x_{k - 1}) = x_{k} \Big\}. \] 

\noindent 
For any $m \le n < \infty$, we define the following projection operators on the set $\mathscr{R}_{n} (x_{1})$ as follows: 
\begin{displaymath} 
\begin{array}{c c c c l c r c l} 
\Pi_{m} & : & \mathscr{R}_{n} (x_{1}) & \longrightarrow & \mathscr{R}_{m} (x_{1}) & \text{defined by} & \Pi_{m} \left( ( x_{1},\, x_{2},\, \cdots,\, x_{n} ) \right) & = & \left( x_{1},\, x_{2},\, \cdots,\, x_{m} \right); \\ 
\pi_{m} & : & \mathscr{R}_{n} (x_{1}) & \longrightarrow & I & \text{defined by} & \pi_{m} \left( ( x_{1},\, x_{2},\, \cdots,\, x_{n} ) \right) & = & x_{m}. 
\end{array}  
\end{displaymath} 
Allowing a slight abuse of notations, for any $m < \infty$, analogous definitions can be written for the projection operators $\Pi_{m}$ and $\pi_{m}$ defined on $\mathscr{R} (x_{1})$. 
\begin{displaymath} 
\begin{array}{c c c c l c r c l} 
\Pi_{m} & : & \mathscr{R} (x_{1}) & \longrightarrow & \mathscr{R}_{m} (x_{1}) & \text{defined by} & \Pi_{m} \left( ( x_{1},\, x_{2},\, \cdots ) \right) & = & \left( x_{1},\, x_{2},\, \cdots,\, x_{m} \right); \\ 
\pi_{m} & : & \mathscr{R} (x_{1}) & \longrightarrow & I & \text{defined by} & \pi_{m} \left( ( x_{1},\, x_{2},\, \cdots ) \right) & = & x_{m}. 
\end{array}  
\end{displaymath} 

\noindent 
We say $x_{1} \in I$ is a \emph{periodic point} of period $p$ with \emph{periodic orbit} $(x_{1},\, x_{2},\, \cdots,\, x_{p}) \in \mathscr{R}_{p} (x_{1})$ pertaining to the combinatorial data given by some $p$-lettered word $w = (w_{1}\, w_{2}\, \cdots\, w_{p})$ on $N$ letters (the length of $w$ denoted by $|w| = p$), if 
\begin{enumerate} 
\item $\pi_{1} \left( ( x_{1},\, x_{2},\, \cdots ) \right)\ \ =\ \ \pi_{p + 1} \left( ( x_{1},\, x_{2},\, \cdots ) \right)$;   
\item $p$ is the least such positive integer for which the first condition is true, {\it i.e.}, \\ $\pi_{1} \left( ( x_{1},\, x_{2},\, \cdots ) \right)\ \ \ne\ \ \pi_{q} \left( ( x_{1},\, x_{2},\, \cdots ) \right)\ \forall q \le p$; and 
\item there does not exist any distinct $1 \le q, r \le p$ for which \\ $\pi_{q} \left( ( x_{1},\, x_{2},\, \cdots,\, x_{p} ) \right)\ \ =\ \ \pi_{r} \left( ( x_{1},\, x_{2},\, \cdots,\, x_{p} ) \right)$. 
\end{enumerate} 
We identify such a periodic point $x_{1} \in I$ with period $p$ and periodic orbit $(x_{1},\, x_{2},\, \cdots,\, x_{p})$ by looking for periodic blocks of points in $\mathscr{R} (x_{1})$ that satisfy, 
\[ \Pi_{p} \left( ( x_{1},\, x_{2},\, \cdots ) \right)\ \ =\ \ \Pi_{p} \left( ( x_{mp\, +\, 1},\, x_{mp\, +\, 2},\, \cdots ) \right)\ \ \ \forall m \in \mathbb{Z}_{+}. \]

\noindent 
It is a simple observation that corresponding to any $p$-periodic point $x \in I$, there exists a $p$-lettered word $w = (w_{1}\, w_{2}\, \cdots\, w_{p})$ on $N$ letters such that 
\[ T_{w} x\ \ :=\ \ \big( T_{w_{p}} \circ \cdots \circ T_{w_{1}} \big) x\ \ =\ \ x. \] 
For any $n$-lettered word $w = ( w_{1}\, w_{2}\, \cdots\, w_{n} )$, we collect the points satisfying $T_{w} x = x$ in the set ${\rm Fix} (T_{w})$. 
\medskip 

\noindent 
For any $T$-invariant probability measure $\mu$ supported on $X$, let $F, G$ be any two complex-valued integrable functions defined on $X$, the appropriate space denoted by $L^{1} (\mu)$. We say the functions $F$ and $G$ are \emph{cohomologous} to each other if there exists a function $H \in L^{1} (\mu)$ such that $F\, -\, G\ =\ \mathscr{Q} (H)\, -\, H$. If $F$ is cohomologous to the constant function $\mathbf{0}$, then $F$ is called a \emph{coboundary}. For any $F \in \mathcal{C} (X, \mathbb{C})$, we denote and define its $n$-th ergodic sum by 
\begin{equation} 
\label{nthergodicsumF} 
F^{n}\ \ :=\ \ F + \mathscr{Q} (F) + \mathscr{Q}^{2} (F) + \cdots + \mathscr{Q}^{n - 1} (F). 
\end{equation} 
Hence, for any two cohomologous functions $F$ and $G$, it is obvious that their $n$-th ergodic sums evaluated at a periodic point of period $n$ must be the same, {\it i.e.}, 
\[ F^{n}((w,\, x))\ \ =\ \ G^{n}((w,\, x))\ \ \forall (w,\, x) \in {\rm Fix}_{n} (T). \] 

\noindent 
Let $F,\, G \in \mathscr{F}_{\alpha}(X, \mathbb{R})$ with $F$ not being cohomologous to any constant function. Then, by a result due to Ruelle in \cite{dr:78}, the function $t \longmapsto \mathfrak{P}(F + t G)$, where $t \in \mathbb{R}$, is convex and real analytic. Further, from \cite{cp:90} and \cite{pp:90}, we have 
\begin{eqnarray} 
\label{firstderivative} 
\left. \frac{d}{dt} \Big( \mathfrak{P}(F + t G) \Big) \right|_{t\, =\, 0} & = & \int\! G\, d \mu_{F} \\ 
\label{secondderivative} 
\left. \frac{d^{2}}{d t^{2}} \Big( \mathfrak{P}(F + t G) \Big) \right|_{t\, =\, 0} & = & \lim_{n\, \rightarrow\, \infty} \frac{1}{n} \int\! \left( G^{n} - n \int\! G\, d \mu_{F} \right)^{2}\, d \mu_{F}\ >\ 0.
\end{eqnarray}

\noindent 
For simultaneous action of $N$ interval maps at a point $x \in I$, the $n$-th order ergodic sum of any $f \in \mathcal{C} (I, \mathbb{C})$ must be calculated over its appropriate orbit, {\it i.e.}, given a $n$-lettered word $w = (w_{1}\, w_{2}\, \cdots\, w_{n})$ on $N$ letters, we define the composition operator $\mathscr{O}_{w}$ by 
\[ \mathscr{O}_{w} (f)\ \ :=\ \ f \circ T_{w},\ \ \ \ \text{where}\ \  T_{w}\ :=\ T_{w_{n}} \circ \cdots \circ T_{w_{1}}. \] 
Then, the $n$-th order ergodic sum of the function $f$ with respect to the given $n$-lettered word $w$, or (by a slight abuse of notations) the $n$-th order ergodic sum of the function $f$ with respect to any given infinite-lettered word $w = (w_{1}\, w_{2}\, \cdots\, w_{n}\, \cdots)$ on $N$ letters that agree with our $n$-lettered word on the initial $n$ positions is given by 
\begin{eqnarray*} 
f^{n}_{w} (x) & := & \left( f + \mathscr{O}_{(w_{1})} (f) + \cdots + \mathscr{O}_{(w_{1}\, w_{2}\, \cdots\, w_{n - 1})} (f) \right) (x) \\
& = & \left( f + f ( T_{w_{1}} ) + \cdots + f \left( T_{w_{n - 1}} \circ \cdots \circ T_{w_{2}} \circ T_{w_{1}} \right) \right) (x). 
\end{eqnarray*} 

\noindent 
For a given $n$-lettered word $w = (w_{1}\, w_{2}\, \cdots\, w_{n})$ on $N$ letters, we say that two Lebesgue integrable functions $f, g \in L^{1} (\lambda)$ defined on $I$ are \emph{$w$-cohomologous} to each other if there exists an integrable function $h \in L^{1} (\lambda)$, also defined on $I$ such that $f\, -\, g\ =\ \mathscr{O}_{w} h\, -\, h$. Hence, for any two $w$-cohomologous functions $f$ and $g$, we observe that the values of the function evaluated at a periodic point $x_{1}$ of period $n$ with periodic orbit $(x_{1},\, x_{2},\, \cdots,\, x_{n})$ pertaining to the combinatorial data given by the $n$-lettered word $w$, necessarily agree; and so do their $n$-th order ergodic sums. Further, if $f$ is $w$-cohomologous to the constant function $\mathbf{0}$, then $f$ is called a \emph{$w$-coboundary}. 
\bigskip 

\section{Statements of results}
\label{main} 

\noindent 
In this section, we state the main theorems of this paper. The first five results concern the setting of the dynamics of the skew-product map $T$ defined on $X = \Sigma_{N}^{+} \times I$, while the next five results concern the setting of simultaneous dynamics of the concerned interval maps on $I$; thus generalising the situation to maps that evolve with multiple branches. 
\bigskip 

\noindent 
\begin{theorem}[Ergodicity] 
\label{erg1}
The action of $T$ on the product space $X$ is necessarily ergodic with respect to the product measure $\mu$. In other words, the measure of any subset $B \subseteq X$, in the product sigma-algebra of $X$ that satisfies $T^{-1} B = B$, is necessarily $0$ or $1$. 
\end{theorem}
\bigskip 

\noindent 
\begin{theorem}[Rates of recurrence] 
\label{ror1}
Consider $F \in \mathscr{F}_{\alpha}(X, \mathbb{R})$ that satisfies the approximability condition, {\it i.e.}, there exists distinct points $(w_{1},\, x_{1}) \in {\rm Fix}_{p_{1}} (T),\ (w_{2},\, x_{2}) \in {\rm Fix}_{p_{2}} (T)$ and $(w_{3},\, x_{3}) \in {\rm Fix}_{p_{3}} (T)$ with $p_{i} \neq p_{j}$ for $i \neq j$ for which 
\begin{equation} 
\label{dioF}
\frac{F^{p_{2}}((w_{2},\, x_{2}))\; -\; F^{p_{1}}((w_{1},\, x_{1}))}{F^{p_{3}}((w_{3},\, x_{3}))\; -\; F^{p_{1}}((w_{1},\, x_{1}))}\ \ =:\ \ \mathfrak{d}_{1} 
\end{equation} 
is a Diophantine number, {\it i.e.}, there exists $l > 2$ and $m > 0$ such that we have 
\begin{equation} 
\label{dioph} 
\left| \mathfrak{d}_{1} - \frac{p}{q} \right|\ \ \geq\ \ \frac{m}{q^{l}},\ \ \forall p, q \in \mathbb{Z}_{+}. 
\end{equation} 
Further, suppose that there exists a unique real number $\kappa$ such that 
\[ \int\! F d \mu_{\kappa F}\ \ =\ \ 0. \] 
Then, for every $n \in \mathbb{Z}_{+},\ a, b \in \mathbb{R}$ with $a < b$, there exists a positive real constant $C_{1} > 0$ such that  
\begin{equation} 
\label{ror1eq} 
\# \Big\{ (w,\, x) \in {\rm Fix}_{n} (T) : a \leq F^{n}((w,\, x)) \leq b \Big\}\ \ \sim\ \ C_{1}\ \frac{e^{n \mathfrak{P} (\kappa F)}}{\sqrt{n}}\ \int_{a}^{b}\!\! e^{- \kappa t}\, dt. 
\end{equation} 
\end{theorem}
\medskip 

\noindent 
Here, in equation \eqref{ror1eq}, by $a_{n} \sim b_{n}$, we mean that $\displaystyle{\lim\limits_{n\, \to\, \infty} \frac{a_{n}}{b_{n}} = 1}$. 
\bigskip 

\noindent 
\begin{theorem}[Exponential decay of correlations]
\label{doc1} 
For any $F \in \mathscr{F}_{\alpha} (X, \mathbb{R})$ whose equilibrium measure is denoted by $\mu_{F}$, there exists a constant $\vartheta \in (0, 1)$ such that for all $G,\ H \in \mathscr{F}_{\alpha} (X, \mathbb{R})$, we have $C_{2} > 0$ (depending on $G$ and $H$) that satisfies 
\begin{equation} 
\label{doce1}
\left\vert \int\! \mathscr{Q}^{n} (G) H\, d \mu_{F}\; -\; \int\! G\, d \mu_{F}\, \int\! H\, d \mu_{F} \right\vert\ \ \le\ \ C_{2} \vartheta^{n};\ \ \forall n \geq 1. 
\end{equation}
\end{theorem}
\medskip 

\noindent 
The preceding theorem defines the exponential decay of correlations of the distributions $\mathscr{Q}^{n} (G)$ and $H$ with respect to the measure $\mu_{F}$ as $n \rightarrow \infty$. The exponential nature of the decay is evident in the statement of the theorem. The next theorem relates the ergodic sum of $G \in \mathscr{F}_{\alpha} (X, \mathbb{R})$ to what is known as the Brownian motion on some richer probability space. 
\bigskip 

\noindent 
\begin{theorem}[Almost sure invariance principle] 
\label{asip1}
For any $F \in \mathscr{F}_{\alpha} (X, \mathbb{R})$ whose equilibrium measure is denoted by $\mu_{F}$, consider $G \in \mathscr{F}_{\alpha}(X, \mathbb{R})$ satisfying \[ \int G\, d \mu_{F}\ \ =\ \ 0. \] 
Suppose the variance of $H \in \mathscr{F}_{\alpha} (X, \mathbb{R})$ is defined as 
\[ \left( \varsigma (H) \right)^{2}\ \ :=\ \ \lim_{n \rightarrow \infty} \frac{1}{n} \int\! \left( H^{n}\; -\; n\, \int\! H\, d \mu_{F} \right)^{2}\, d \mu_{F}. \] 
Then, there exists a H\"{o}lder continuous function $\Phi \in \mathscr{F}_{\alpha}(X, \mathbb{R})$ cohomologous to $G$, a one-dimensional Brownian motion $\Big\{ \mathfrak{B}(t) \Big\}_{t\, \ge\, 0}$ with variance $t \left( \varsigma (G) \right)^{2}$ and a sequence of random variables $\big\{ \mathfrak{Y}_{n} : \Omega \longrightarrow \mathbb{R} \big\}_{n\, \ge\, 0}$ such that $\big\{ \mathfrak{Y}_{n} \big\}_{n\, \ge\, 0}$ and $\big\{ \Phi^{n} \big\}_{n\, \ge\, 0}$ are equal in distribution and given any $\delta > 0$, 
\[ \mathfrak{Y}_{\lfloor t \rfloor}(\omega)\ \ =\ \ \mathfrak{B} (t) (\omega)\; +\; O(t^{\frac{1}{4}\, +\, \delta}),\ \ \forall t \geq 0,\ \ \mu_{F} \text{-a.e.}, \]
provided $G$ is not a coboundary.
\end{theorem}
\medskip 

\noindent 
The almost sure invariance principle leads to a few important corollaries such as the central limit theorem and the law of iterated logarithms. 
\bigskip 

\noindent 
\begin{theorem}[Central limit theorem and law of iterated logarithms] 
\label{osp1}
For any $F \in \mathscr{F}_{\alpha} (X, \mathbb{R})$ whose equilibrium measure is denoted by $\mu_{F}$, consider $G \in \mathscr{F}_{\alpha}(X, \mathbb{R})$ satisfying 
\[ \int\! G\, d \mu_{F}\ \ =\ \ 0. \] 
Suppose $G$ is not a coboundary. Then, 
\begin{enumerate} 
\item $G$ satisfies the central limit theorem, {\it i.e.}, $\frac{1}{\sqrt{n}} G^{n}$ converges in distribution to a normal distribution with mean zero and variance $\left( \varsigma (G) \right)^{2}$ as $n \rightarrow \infty$. 
\item $G$ satisfies the law of iterated logarithms, {\it i.e.}, 
\[ \limsup_{n\, \rightarrow\, \infty} \frac{G^{n} ((w,\, x))}{\varsigma (G) \sqrt{2n \log \log n}}\ \ =\ \ 1\ \ \mu_{F}\text{-a.e.} \] 
\end{enumerate}
\end{theorem}
\medskip 

\noindent 
In the next five theorems, we state theorems captioned under the same titles, however, by suppressing the first co-ordinate of $X$ and looking at a genuine simultaneous action of the finitely many interval maps under consideration. 
\bigskip 

\noindent 
\begin{theorem}[Ergodicity] 
\label{erg2}
Consider the interval maps $T_{1},\, T_{2},\, \cdots,\, T_{N}$ that act simultaneously on the interval $I$. Let $\lambda$ denote the Lebesgue measure on $\mathbb{R}$. Then, for any real-valued Lebesgue integrable function $f \in L^{1} (\lambda)$, for $\lambda$-a.e. $x \in I$, we have 
\begin{equation} 
\lim_{n\, \to\, \infty}\, \frac{1}{n}\, \frac{1}{N^{n}}\, \sum_{w\; :\; |w|\, =\, n} \Big[ f\, +\, \mathscr{O}_{(w_{1})} (f)\, +\, \cdots\, +\, \mathscr{O}_{(w_{1}\, w_{2}\, \cdots\, w_{n - 1})} (f) \Big] (x)\ \ =\ \ \int_{0}^{1}\! f\, d \lambda. 
\end{equation} 
\end{theorem}
\bigskip 

\noindent 
\begin{theorem}[Rates of recurrence] 
\label{ror2}
Consider $f \in \mathscr{F}_{\alpha}(I, \mathbb{R})$ that satisfies the approximability condition, {\it i.e.}, there exists distinct periodic points $x,\ y$ and $z$ in $I$ with distinct periods $p_{x},\ p_{y}$ and $p_{z}$, pertaining to the combinatorial data given by $w_{x},\ w_{y}$ and $w_{z}$ such that 
\begin{equation} 
\label{dioph2} 
\frac{f^{p_{y}}_{w_{y}} (y)\ -\ f^{p_{x}}_{w_{x}} (x)}{f^{p_{z}}_{w_{z}} (z)\ -\ f^{p_{x}}_{w_{x}} (x)}\ \ =\ \ \mathfrak{d}_{2} 
\end{equation} 
is a Diophantine number. Further, suppose there exists unique $\kappa > 0$ such that 
\[ \int\! f d\mathfrak{m}_{\kappa f}\ \ =\ \ 0. \] 
Then, for every $n \in \mathbb{Z}_{+},\ a, b \in \mathbb{R}$ with $a < b$, there exists a positive real constant $C_{3} > 0$ such that 
\[ \sum\limits_{w\; :\; |w|\, =\, n} \# \big\{ x \in {\rm Fix} (T_{w})\ :\ a \leq f^{n}_{w}(x) \leq b \big\}\ \ \sim\ \ C_{3} \frac{e^{n \mathbb{P}(\kappa f)}}{\sqrt{n}}\ \int_{a}^{b}\! e^{-\kappa t}\, d t. \]
\end{theorem}
\bigskip 

\noindent 
\begin{theorem}[Exponential decay of correlations]
\label{doc2}
Let $\lambda$ denote the Lebesgue measure on $I$. Then, there exist a constant $\vartheta \in (0, 1)$ such that for all $\alpha$-H\"older continuous functions $g,\ h \in \mathscr{F}_{\alpha}(I, \mathbb{R})$, we have $C_{4} > 0$ (depending on $g,\ h$ and some $n$-lettered word $w$) that satisfies
\begin{equation}
\label{doce2}
\left\vert \int\! \mathscr{O}_{w} (g) h\, d \lambda\; -\; \int\! g\, d \lambda \, \int\! h\, d \lambda  \right\vert\ \ \le\ \ C_{4} \vartheta^{n};\ \ \forall n \geq 1.
\end{equation}
\end{theorem}
\bigskip 

\noindent 
\begin{theorem}[Almost sure invariance principle] 
\label{asip2}
Let $\lambda$ denote the Lebesgue measure on $I$. For any $g \in \mathscr{F}_{\alpha} (I, \mathbb{R})$ with 
\[ \int\! g\, d \lambda\ \ =\ \ 0, \] 
and $w = (w_{1}\, w_{2}\, \cdots) \in \Sigma_{N}^{+} $, assume that the variance of $g$ with respect to the word $w$ denoted by $\big( \varsigma_{w} (g) \big)^{2}$ and defined by 
\[ \big( \varsigma_{w} (g) \big)^{2}\ \ :=\ \ \lim\limits_{n\, \to\, \infty} \frac{1}{n} \int\! \left( g_{w}^{n} \right)^{2}\, d \lambda\ \ >\ \ 0. \] 
Then, there exists a probability space $( \Omega,\; \mathscr{A},\; \nu )$, a sequence of random variables $\big\{ Y_{w}^{n} \big\}_{n\, \ge\, 0}$ and a standard Brownian motion $\big\{ \mathfrak{B}^{*} (t) \big\}_{t\, \ge\, 0}$ such that $g_{w}^{n}$ and $Y_{w}^{n}$ are equal in distribution and given any $\delta > 0$, 
\begin{equation} 
\label{varsigmawng} 
Y_{w}^{n} (\omega)\ -\ \mathfrak{B}^{*} \left( \left( \varsigma_{w}^{(n)} (g) \right)^{2} \right) (\omega)\ \ =\ \ O \big( n^{\frac{1}{4}\, +\, \delta} \big),\ \ \ \ \nu\text{-a.e.}\ \ \ \text{where}\ \ \left( \varsigma_{w}^{(n)} (g) \right)^{2}\ \ =\ \ \int\! \left( g_{w}^{n} \right)^{2}\, d \lambda, 
\end{equation} 
provided $g_{\Pi_{n} (w)}$ is not a $\Pi_{n} (w)$-coboundary for any $n \ge 1$. 

\end{theorem}
\bigskip 

\noindent 
\begin{theorem}[Central limit theorem and law of iterated logarithms] 
\label{osp2}
For a function $g \in \mathscr{F}_{\alpha} (I, \mathbb{R})$ that satisfies
\[ \int\! g\, d \lambda\ \ =\ \ 0, \] 
suppose $w \in \Sigma_{N}^{+} $ is such that the variance $\big( \varsigma_{w} (g) \big)^{2} > 0$ and $g_{\Pi_{n} (w)}$ is not a $\Pi_{n} (w)$-coboundary for any $n \ge 1$. Then,
\begin{enumerate}
\item $g$ satisfies the central limit theorem \textit{i.e.,} $\frac{1}{\sqrt{n}} g_{w}^{n} $ converges in distribution to a normal distribution with mean zero and variance $\big( \varsigma_{w} (g) \big)^{2}$ as $n \to \infty$. 
\item $g$ satisfies the law of iterated logarithms, {\it i.e.}, 
\[ \limsup_{n\, \rightarrow\, \infty} \frac{g_{w}^{n} (x)}{\varsigma_{w} (g) \sqrt{2n \log \log n}}\ \ =\ \ 1\ \ \lambda\text{-a.e.} \] 
\end{enumerate}
\end{theorem} 
\bigskip 

\section{A book-keeping mechanism, $\Sigma_{N}^{+}$} 
\label{SigmaNplus}

In this section, we define the space $\Sigma_{N}^{+}$ and discuss certain properties that will be useful in the sequel. Interested readers may refer to \cite{bk:98}, for more details on this space. Recall the definition of $\Sigma _{N}^{+}$ from section \eqref{prelims}, 
\[ \Sigma_{N}^{+}\ \ :=\ \ S^{\mathbb{Z}_{+}}\ \ =\ \ \big\{ 1,\, 2,\, \cdots,\, N \big\}^{\mathbb{Z}_{+}}\ \ =\ \ \Big\{ w = \left( w_{1}\, w_{2}\, \cdots\, w_{n}\, \cdots \right)\ :\ w_{i} \in \left\{ 1,\, 2,\, \cdots,\, N \right\} \Big\}. \] 
Observe that one can define the maps $\Pi_{m}$ and $\pi_{m}$ on the symbolic space $S^{n}$ as well as $\Sigma_{N}^{+}$, analogous to its definitions on $I^{n}$ and $I^{\mathbb{Z}_{+}}$. We make use of the same to define a metric between any words $v, w \in \Sigma_{N}^{+}$. Fix any $\theta \in (0, 1)$, and define 
\[ d_{\Sigma_{N}^{+}} (v,\, w)\ \ :=\ \ \theta^{n (v,\, w)},\ \ \text{where}\ \ n (v,\, w)\ :=\ \sup \Big\{ k \in \mathbb{Z}_{+}\ :\ \Pi_{k} (v)\ =\ \Pi_{k} (w) \Big\}. \] 
Here, we define $n ( v, v ) := \infty$, thereby $d_{\Sigma_{N}^{+}} ( v, v ) = 0$. Thus, it is clear that we have a family of metrics on the space $\Sigma_{N}^{+}$. The discrete topology that separates any two distinct symbols on the set $\{ 1, 2, \cdots, N \}$ accords a product topology on $\Sigma_{N}^{+}$ with which the above described family of metrics is compatible. We shall fix a value of $\theta$, according to our need in a later section. In this topology, the cylinder sets given by fixing a finite set of co-ordinates, are the sets that are both closed and open. For ease of explanations, we shall always consider cylinder sets whose co-ordinates are fixed from the first co-ordinate onwards, for example, a cylinder set of length $m$ looks like 
\[ \big[ v_{1}\, v_{2}\, \cdots\, v_{m} \big]\ \ =\ \ \Big\{ w \in \Sigma_{N}^{+}\ :\ \Pi_{m} (w)\ =\ (v_{1}\, v_{2}\, \cdots\, v_{m}) \Big\}. \] 
These cylinder sets form a basis for the $\sigma$-algebra on $\Sigma_{N}^{+}$ on which one could define a measure for $\Sigma_{N}^{+}$. An easily describable measure on the space, $\Sigma_{N}^{+}$ is the Bernoulli measure defined thus. For any fixed probability vector $p = \big( p_{1},\, p_{2},\, \cdots,\, p_{N} \big)$, the measure is defined as 
\[ \mu \Big( \big[ v_{1}\, v_{2}\, \cdots\, v_{m} \big] \Big)\ \ =\ \ p_{v_{1}}\; p_{v_{2}}\; \cdots\; p_{v_{m}}. \] 
Observe that the shift map $\sigma$ defined on $\Sigma_{N}^{+}$ satisfies the properties asked for with respect to the topology defined on $\Sigma_{N}^{+}$. Further, $\Sigma_{N}^{+}$ is a compact metric space with topological dimension $0$. 
\medskip 

\noindent 
We now appeal to Tychonoff and accord some structure on $X$. For any two points $(w,\, x)$ and $(v,\, y)$ in $X$, we define a metric 
\[ d_{X} ((w,\, x),\, (v,\, y))\ \ :=\ \ \max \Big\{ d_{\Sigma_{N}^{+}} (w,\, v),\; |x - y| \Big\}. \] 
Thus, we work with the appropriate product topology and the product $\sigma$-algebra and the product measure, while we work with $X$. 

\section{Various Ruelle operators} 
\label{vro}

Making use of the Ruelle operator, as given in \cite{pp:90} for every $1 \le d \le N$, we define a Ruelle operator for the skew-product map in this section. Later, we consider each of these Ruelle operators to define a collective Ruelle operator for the case of simultaneous action of all these maps. 
\medskip 

\noindent 
For every $1 \le d \le N$, fix $f \in \mathcal{C} (I, \mathbb{C})$, consider $\mathcal{L}_{f}^{(d)} : \mathcal{C} (I, \mathbb{C}) \longrightarrow \mathcal{C} (I, \mathbb{C})$ given by 
\begin{equation} 
\label{ruelled} 
\left( \mathcal{L}_{f}^{(d)} g \right) (x)\ \ :=\ \ \sum_{T_{d} y\, =\, x} e^{f(y)} g(y). 
\end{equation}  
Observe that this definition entails the following iterative formula given by, 
\begin{equation} 
\label{iterue} 
\left( \left( \mathcal{L}_{f}^{(d)} \right)^{\!\circ n} g \right) (x)\ \ :=\ \ \sum_{T_{d}^{n} y\, =\, x} e^{f^{n}_{(d\, d\, \cdots\, d)} (y)} g(y). 
\end{equation} 
This is the usual Ruelle operator, as defined in \cite{pp:90}. For such an operator $\mathcal{L}_{f}^{(d)}$, we have the Ruelle operator theorem, as stated in \cite{pp:90}. 
\medskip 

\noindent 
\begin{theorem}[\cite{pp:90, du:91}] 
\label{rotd}
Suppose $f \in \mathscr{F}_{\alpha} (I, \mathbb{R})$. Then, the Ruelle operator $\mathcal{L}_{f}^{(d)}$ has a simple maximal positive eigenvalue, $\rho^{(d)}$. The remainder of the spectrum lies in a disc of radius strictly smaller than $\rho^{(d)}$. The eigenfunction $\phi_{d}$ corresponding to the maximal eigenvalue is strictly positive. Further, there exists an eigenmeasure corresponding to the maximal eigenvalue, in the space of all $T_{d}$-invariant probability measures supported on $I$, for the dual operator $\left( \mathcal{L}_{f}^{(d)} \right)^{\!*}$. 
\end{theorem} 
\medskip 

\noindent 
Appealing to variational principles studied by several authors including Bowen in \cite{rb:73}, Ruelle in \cite{dr:78} and Parry and Pollicott in \cite{pp:90}, we know that the maximal eigenvalue for the Ruelle operator $\mathcal{L}_{f}^{(d)}$, for $f \in \mathscr{F}_{\alpha} (I, \mathbb{R})$ can also be described using the $d$-pressure function, {\it i.e.}, $\rho^{(d)} = e^{\mathcal{P}^{(d)} (f)}$, where 
\[ \mathcal{P}^{(d)} (f)\ \ :=\ \ \sup \left\{ h_{m^{(d)}} (T_{d}) + \int\! f\, d m^{(d)} \right\}.  \] 
The supremum in the above definition is taken over all $T_{d}$-invariant probability measures supported on $I$ and $h_{m^{(d)}} (T_{d})$ is the measure theoretic entropy of $T_{d}$ with respect to the measure $m^{(d)}$. Then, the existence of the unique equilibrium measure, denoted by $m_{f}^{(d)}$, that realises the supremum in the definition of $d$-pressure is assured by Denker and Urbanskii in \cite{du:91}. Further, the variational principle states that this unique equilibrium measure, $m_{f}^{(d)}$ is equivalent to the eigenmeasure corresponding to the maximal eigenvalue $\rho^{(d)}$ for the dual operator $\left( \mathcal{L}_{f}^{(d)} \right)^{\!*}$, as given in theorem \eqref{rotd}. 
\medskip 

\noindent 
Taking cue from definition \eqref{ruelled}, we define the Ruelle operator for the skew-product setting as follows: Fix $F \in \mathcal{C} (X, \mathbb{C})$ and consider $\mathfrak{L}_{F} : \mathcal{C} (X, \mathbb{C}) \longrightarrow \mathcal{C} (X, \mathbb{C})$ given by 
\begin{equation} 
\label{ruelleT} 
\left( \mathfrak{L}_{F} G \right) ((w,\, x))\ \ :=\ \ \sum_{T ((v,\, y))\, =\, (w,\, x)} e^{F((v,\, y))} G((v,\, y)). 
\end{equation} 
It is a simple observation that the iterates of the Ruelle operator $\mathfrak{L}_{F}$ agrees with the appropriate iterative formula given in \eqref{iterue}. We merely state the same here. 
\[ \left( \left( \mathfrak{L}_{F} \right)^{\!\circ n} G \right) ((w,\, x))\ \ :=\ \ \sum_{T^{n} ((v,\, y))\, =\, (w,\, x)} e^{F^{n} ((v,\, y))} G((v,\, y)). \] 
Further, the operator $\mathfrak{L}_{F}$ satisfies the properties mentioned in the Ruelle operator theorem, as mentioned in theorem \eqref{rotd}, whenever $F \in \mathscr{F}_{\alpha} (X, \mathbb{R})$. We include the statement of the theorem, in the context of the skew-product map, for readers' convenience. 
\medskip 

\noindent 
\begin{theorem} 
\label{rotskew}
Suppose $F \in \mathscr{F}_{\alpha} (X, \mathbb{R})$. Then, the Ruelle operator $\mathfrak{L}_{F}$ has a simple maximal eigenvalue at $\varrho = e^{\mathfrak{P} (F)}$. The remainder of the spectrum lies in a disc of radius strictly smaller than $e^{\mathfrak{P} (F)}$. The eigenfunction $\Phi$ corresponding to the maximal eigenvalue is strictly positive. Further, the eigenmeasure corresponding to the maximal eigenvalue for the dual operator $\left( \mathfrak{L}_{F} \right)^{\!*}$ is equivalent to the equilibrium measure $\mu_{F}$, that realises the supremum in the definition of pressure, as stated in equation \eqref{pressure}. 
\end{theorem} 
\medskip 

\noindent 
Fixing $f \in \mathcal{C} (I, \mathbb{C})$, we now define a third Ruelle operator, denoted by 
\[ \mathbb{L}_{f}\ :\ \mathcal{C} (I, \mathbb{C})\ \longrightarrow\ \mathcal{C} (I, \mathbb{C}), \] 
that captures the idea of simultaneous action of interval maps. Making sense of the $n$-th iterate of the Ruelle operators $\mathcal{L}_{f}^{(d)}$ and $\mathfrak{L}_{F}$ that captures the set of all $n$-th order pre-images of the point where the operator acts and taking the $n$-th ergodic sum over each of those orbits, we are inclined to define the $n$-th iterate of the Ruelle operator $\mathbb{L}_{f}$ acting on a point $x$ by considering those points that would reach $x$ in $n$ steps, by the action of a combination of $n$ many maps from the collection $\{ T_{1},\, T_{2},\, \cdots,\, T_{N} \}$ and taking the $n$-th ergodic sum dictated by all such $n$-lettered words. The same can be expressed as 
\begin{eqnarray*} 
\left( \left( \mathbb{L}_{f} \right)^{\!\circ n} g \right) (x) & = & \sum_{w\, :\, |w|\, =\, n} \sum_{\left( T_{w_{n}} \circ T_{w_{n - 1}} \circ \cdots \circ T_{w_{1}} \right) y\, =\, x} e^{f^{n}_{w} (y)} g(y) \\
& = & \sum_{w_{n}\, =\, 1}^{N} \cdots \sum_{w_{1}\, =\, 1}^{N} \sum_{\left( T_{w_{n}} \circ T_{w_{n - 1}} \circ \cdots \circ T_{w_{1}} \right) y\, =\, x} e^{f (y)\, +\, f (T_{w_{1}} y)\, +\, \cdots\, +\, f ( T_{w_{n - 1}} \circ \cdots \circ T_{w_{1}} y)} g(y). 
\end{eqnarray*} 

\noindent 
This understanding paves the way for us to define the Ruelle operator, in this case as 
\begin{equation} 
\label{ruellesimul} 
\mathbb{L}_{f} g (x)\ \ :=\ \ \sum_{d\, =\, 1}^{N} \sum_{T_{d} y\, =\, x} e^{f(y)} g(y)\ \ =\ \ \sum_{d\, =\, 1}^{N} \mathcal{L}_{f}^{(d)} g (x). 
\end{equation} 

\section{Spectrum of the operators $\mathfrak{L}_{F}$ and $\mathbb{L}_{f}$} 
\label{spectrum} 

\noindent 
In this section, we establish a relationship between the Ruelle operators $\mathfrak{L}_{F}$ and $\mathbb{L}_{f}$, as defined in equations \eqref{ruelleT} and \eqref{ruellesimul}. 
\medskip 

\noindent 
Let 
\begin{equation} 
\label{qoff} 
Q\ :\ \mathscr{F}_{\alpha} (I, \mathbb{C}) \longrightarrow \mathscr{F}_{\alpha} (X, \mathbb{C})\ \ \text{be defined as}\ \ \left( Q(f) \right) (w,\, x)\ \ :=\ \ f (x). 
\end{equation} 
Then, for any H\"{o}lder continuous function $f \in \mathscr{F}_\alpha (I, \mathbb{C})$, 
\begin{eqnarray*} 
\big| Q(f) ((w,\, x))\ -\ Q(f)((v,\, y)) \big| & = & \big| f(x) - f(y) \big| \\ 
& \le & M_{f} \left| x - y \right|^{\alpha} \\ 
& \le & M_{f}\, \big[ d_{X} ((w,\, x),\ (v,\, y)) \big]^{\alpha}, 
\end{eqnarray*}
for some $M_{f} > 0$ and for any $0 < \theta < 1$, on which the metric on $\Sigma_{N}^{+}$ depends. Thus, the map $Q$ is well defined. Further, the above inequality also proves that the H\"{o}lder constant $M_{f}$ remains unperturbed for the function $Q(f)$ in the product space as well, {\it i.e.}, $M_{f} \equiv M_{Q(f)}$. Moreover, it is clear from the definition of the various Ruelle operators that 
\[ Q \left( \mathbb{L}_{f} g \right)\ \ =\ \ \mathfrak{L}_{Q(f)} Q(g),\ \ \forall f,\, g \in \mathscr{F}_{\alpha} (I, \mathbb{C}). \] 
Thus, the action of $\mathbb{L}_{f}$ is similar to that of the action of $\mathfrak{L}_{Q(f)}$ restricted to the subspace, ${\rm Image} (Q) \subseteq \mathscr{F}_{\alpha} (X, \mathbb{C})$. As a consequence, we can relate the spectrum of $\mathbb{L}_{f}$ and $\mathfrak{L}_{Q(f)}$. The following lemma narrates the same. 
\medskip 

\noindent 
\begin{lemma} 
\label{rotsa} 
For some fixed $f \in \mathscr{F}_{\alpha}(I, \mathbb{C})$, let $\Phi \in \mathscr{F}_{\alpha} (X, \mathbb{C})$ be an eigenfunction of $\mathfrak{L}_{Q(f)}$ with corresponding eigenvalue $\varrho$, {\it i.e.}, $\mathfrak{L}_{Q(f)} \Phi = \varrho\, \Phi$. Then, $\Phi \in {\rm Image}(Q) \subset \mathscr{F}_{\alpha} (X, \mathbb{C})$. 
\end{lemma}

\noindent 
\begin{prooof} 
In order to prove this lemma, it is sufficient to prove that the eigenfunction $\Phi$ is independent of the first co-ordinate, {\it i.e.}, $\Phi ((v,\, x)) = \Phi ((w,\, x)),\ \forall v,\, w \in \Sigma_{N}^{+}$ and $x \in [0,\, 1]$. In fact, we prove that given any $\epsilon > 0$, there exists $M_{\epsilon} \in \mathbb{Z}_{+}$ such that 
\[ \Big\vert \left( \mathfrak{L}_{Q(f)} \right)^{\!\circ n} \Phi ((v,\, x))\ -\ \left( \mathfrak{L}_{Q(f)} \right)^{\!\circ n} \Phi ((w,\, x)) \Big\vert\ \ \le\ \ \varrho^{n} \epsilon,\ \ \forall n \ge M_{\epsilon}. \] 
A mere application of the eigenfunction equation to the above inequality, then yields 
\[ \varrho^{n} \Big| \Phi ((v,\, x))\ -\ \Phi ((w,\, x)) \Big|\ \ \le\ \ \varrho^{n} \epsilon,\ \ \forall n \ge M_{\epsilon}, \] 
implying 
\[ \Big| \Phi ((v,\, x))\ -\ \Phi ((w,\, x)) \Big|\ \ \le\ \ \epsilon. \] 
The proof is then complete, appealing to the arbitrary choice that we can make for $\epsilon$. 
\medskip 

\noindent 
Using the definition of the Ruelle operator, we have 
\begin{eqnarray*} 
& & \Big| \left( \mathfrak{L}_{Q(f)} \right)^{\!\circ n} \Phi ((v,\, x)) - \left( \mathfrak{L}_{Q(f)} \right)^{\!\circ n} \Phi ((w,\, x)) \Big| \\ 
& = & \left| \sum_{u\; :\; |u|\, =\, n} \left( \sum_{T^{n} ((u v,\, y))\, =\, (v,\, x)} e^{(Q(f))^{n} (u v,\, y)} \Phi ((u v,\, y)) \right. \right. \\ 
& & \hspace{+2in} \left. \left. -\ \sum_{T^{n} ((u w,\, y))\, =\, (w,\, x)} e^{(Q(f))^{n} (u w,\, y)} \Phi ((u w,\, y)) \right) \right| \\ 
& \le & \sum_{u\; :\; |u|\, =\, n} \sum_{T_{u_{n}} \circ \cdots \circ T_{u_{1}} y\, =\, x} \Big| e^{f_{u}^{n} (y)} \Big| \Big| \Phi ((u v,\, y)) - \Phi ((u w,\, y)) \Big| \\ 
& \leq & e^{n \| f \|_{\infty}}\, M_{\Phi}\, N^{n}\, \theta^{n \alpha} \\ 
& = & \varrho^{n}\, M_{\Phi}\, \left( \frac{e^{\| f \|_{\infty}}\, N\, \theta^{\alpha}}{\varrho} \right)^{n}.  
\end{eqnarray*}

\noindent 
Once we choose the functions $f$ and $\Phi$, observe that the quantities $\| f \|_{\infty}$ and $M_{\Phi}$ are determined. Thus, we can only rely on our choice of $\theta$ to make the above estimate, as small as necessary. We therefore fix $\theta$ in such fashion that 
\[ \frac{e^{\| f \|_{\infty}}\, N\, \theta^{\alpha}}{\varrho}\ \ <\ \ 1. \] 
Further, we remark that we shall use the same $\theta$, so fixed to suit our purpose in the above inequality, in the remainder of this paper. Now, it is clear that we have some threshold, say $M_{\epsilon}$ such that 
\[ \Big| \left( \mathfrak{L}_{Q(f)} \right)^{\!\circ n} \Phi ((v,\, x))\ -\ \left( \mathfrak{L}_{Q(f)} \right)^{\!\circ n} \Phi ((w,\, x)) \Big|\ \ \le\ \ \varrho^{n} \epsilon,\ \ \forall n \ge M_{\epsilon}. \] 
\end{prooof}

\noindent 
A careful reader may observe that lemma \eqref{rotsa} is merely the statement of the Ruelle operator theorem stated for the simultaneous dynamics of finitely many interval maps. In particular, $\Phi \in {\rm Image}(Q)$ implies that there exists $\phi \in \mathscr{F}_{\alpha} (I, \mathbb{C})$ such that $Q(\phi) = \Phi$ and $\mathbb{L}_{f} \phi = \varrho \phi$. Thus, the set of eigenvalues of $\mathfrak{L}_{Q(f)}$ and $\mathbb{L}_{f}$ remains equal. In particular, the simple maximal eigenvalue $\varrho$ of $\mathfrak{L}_{Q(f)}$ is also the simple maximal eigenvalue of $\mathbb{L}_{f}$. 

\section{Normalising the operators $\mathcal{L}_{f}^{(d)},\ \mathfrak{L}_{F}$ and $\mathbb{L}_{f}$}
\label{Normalising}

\noindent 
For technical convenience, we normalise the Ruelle operators $\mathcal{L}_{f}^{(d)},\ \mathfrak{L}_{Q(f)}$ and $\mathbb{L}_{f}$, in different ways that would suit our purposes to prove the main theorems. 
\medskip 

\noindent 
For any $f \in \mathscr{F}_{\alpha} (I, \mathbb{R})$, we define the normalised Ruelle operator by considering the Ruelle operator, as defined in \eqref{ruelled}, corresponding to the function 
\[ \widetilde{f}_{d}\ \ :=\ \ f\; +\; \log \phi_{d}\; -\; \log \phi_{d} \circ T_{d}\; -\; \mathcal{P}^{(d)} ( f ), \] 
where $\phi_{d}$ is the eigenfunction corresponding to the maximal eigenvalue of the operator $\mathcal{L}_{f}^{(d)}$. Thus, 
\[ \widetilde{\mathcal{L}}_{f}^{(d)}\, \mathbf{1} (x)\ \ :=\ \ \mathcal{L}_{\widetilde{f}_{d}}^{(d)}\, \mathbf{1} (x)\ \ =\ \ \mathbf{1} (x). \] 
\medskip 

\noindent 
For any $f \in \mathscr{F}_{\alpha} (I, \mathbb{R})$, we define the normalised Ruelle operator in the skew-product setting by considering the Ruelle operator, as defined in \eqref{ruelleT}, corresponding to the function 
\[ \widetilde{Q (f)}\ \ :=\ \ Q(f)\; +\; \log \Phi\; -\; \log  \Phi \circ T\, , \] 
where $\Phi$ is the eigenfunction corresponding to the maximal eigenvalue $\varrho = e^{\mathfrak{P} ( Q(f))}$ of the operator $\mathfrak{L}_{Q(f)}$. Thus, this normalisation effects the eigenfunction of the normalised operator to be equal to the constant function $\mathbf{1}$, however, leaves the the maximal eigenvalue $\varrho$ unchanged, {\it i.e.}, 
\[ \widetilde{\mathfrak{L}}_{Q(f)}\, \mathbf{1} ((w,\, x))\ \ :=\ \ \mathfrak{L}_{\widetilde{Q(f)}}\, \mathbf{1} ((w,\, x))\ \ =\ \ \varrho\, \mathbf{1} ((w,\, x)). \] 
In particular, the operator $\varrho^{-1} \widetilde{\mathfrak{L}}_{Q(f)}$ has $1$ as its maximal eigenvalue with corresponding eigenfunction $\mathbf{1}$. Further, since by definition, $Q(f)$ and $\widetilde{Q(f)}$ are cohomologous to each other, their ergodic sums are preserved. 
\medskip 

\noindent 
Further, once we normalise the operators $\mathcal{L}_{f}^{(d)}$ and $\mathfrak{L}_{Q(f)}$ as prescribed, we also observe that the equilibrium measure with respect to the H\"{o}lder continuous functions $f \in \mathscr{F}_{\alpha} (I, \mathbb{R})$ for the map $T_{d}$ and $Q(f) \in \mathscr{F}_{\alpha} (X, \mathbb{R})$ for the skew-product map $T$, namely $m_{f}^{(d)}$ and $\mu_{Q(f)}$ are nothing but the eigenmeasures corresponding to the maximal eigenvalues of the normalised operators $\widetilde{\mathcal{L}}_{f}^{(d)}$ and $\varrho^{-1} \widetilde{\mathfrak{L}}_{Q(f)}$ respectively. Thus, we use the notations $m_{f}^{(d)}$ and $\mu_{Q(f)}$ for the equilibrium measure, as well as the eigenmeasure corresponding to the maximal eigenvalue of the normalised operators $\widetilde{\mathcal{L}}_{f}^{(d)}$ as well as $\varrho^{-1} \widetilde{\mathfrak{L}}_{Q(f)}$, respectively. 
\medskip 

\noindent 
Moreover, we define the normalised operator in the setting of simultaneous dynamics thus: For any $f \in \mathscr{F}_{\alpha} (I, \mathbb{R})$, let 
\[ \widetilde{\mathbb{L}}_{f} g (x)\ \ :=\ \ \frac{1}{\phi (x)}\, \mathbb{L}_{f + \log \phi} g (x)\ \ =\ \ \frac{1}{\phi (x)}\, \sum_{d\, =\, 1}^{N} \sum_{T_{d} y\, =\, x} e^{f(y)\, +\, \log \phi (y)} g(y), \] 
where $\phi$ is the strictly positive eigenfunction corresponding to the maximal eigenvalue of the operator $\mathbb{L}_{f}$, that we denote by $e^{\mathbb{P} (f)}$, where $\mathbb{P} (f)$ is the \emph{pressure} of the function $f$, in this setting. It must be noted that this normalisation does not change the maximal eigenvalue, but has an effect on the corresponding eigenfunction, making it to be $\mathbf{1}$. Further, this normalisation also takes a toll in this setting; $\widetilde{\mathbb{L}}_{f}$ is no longer a Ruelle operator, but merely a bounded linear operator. 
\medskip 
 
\noindent
The operator $Q : \mathcal{C} (I, \mathbb{C}) \longrightarrow \mathcal{C} (X, \mathbb{C})$, as defined in \eqref{qoff}, has a natural transpose 
\[ Q^{*}\ :\ \mathcal{C}^{*} (X, \mathbb{C}) \longrightarrow \mathcal{C}^{*} (I, \mathbb{C}). \] 
On a restricted space, this transpose gives us the map, 
\[ Q^{*}\ :\ \big\{ \mu\ :\ \mu\ \text{is a probability measure on}\ X \big\} \longrightarrow \big\{ \mathfrak{m}\ :\ \mathfrak{m}\ \text{is a probability measure on}\ I \big\}, \] 
defined by 
\[ \int\! f\, d \left( Q^{*} \mu \right)\ \ :=\ \ \int\! Q(f)\, d \mu,\ \ \forall f \in \mathcal{C} (I, \mathbb{C}). \] 
For $f \in \mathscr{F}_{\alpha}(X, \mathbb{R})$, we define $\mathfrak{m}_{f} := Q^{*}(\mu_{Q(f)})$, thereby, 
\[ \int_{I}\! g\, d \mathfrak{m}_{f}\ \ =\ \ \int_{X}\! Q(g)\, d \mu_{Q(f)},\ \ \forall g \in \mathscr{F}_{\alpha} (I, \mathbb{R}). \]
Thus, it is an easy observation that as a consequence of the definitions of $\mathfrak{m}_{f}$ and $\widetilde{\mathbb{L}}_{f},\ \mathfrak{m}_{f}$ coincides with the eigenmeasure corresponding to maximal eigenvalue of the operator $e^{-\mathbb{P}(f)} \left( \widetilde{\mathbb{L}}_{f} \right)^{*}$. 
\medskip
 
\noindent 
In the sequel, we will be interested in a complex perturbation of the operators $\widetilde{\mathfrak{L}}_{Q(f)}$ and $\widetilde{\mathbb{L}}_{f}$ for $f \in \mathscr{F}_{\alpha} (I, \mathbb{R})$, {\it i.e.}, for $\zeta = \kappa + i \xi \in \mathbb{C}$, we consider the operators $\mathfrak{L}_{\kappa Q(f)}$ and $\mathbb{L}_{\kappa f}$, normalise the same as explained above respectively and then perturb the operator. We denote them respectively by $\widetilde{\mathfrak{L}}_{\zeta Q(f)}$ and $\widetilde{\mathbb{L}}_{\zeta f}$. We conclude this section with a few lemmas that provide a bound on the operator norm for the iterates of the normalised, yet perturbed operators. 
\medskip 

\noindent 
\begin{lemma}[\cite{pp:90}] 
For any $F \in \mathscr{F}_{\alpha} (X, \mathbb{R})$, there exists a positive constant $C_{5} > 0$, such that for every $n \ge 0$, and any $G \in \mathscr{F}_{\alpha} (X, \mathbb{C})$, we have 
\[ \bigg\Vert \left( \varrho^{-1} \widetilde{\mathfrak{L}}_{\zeta F} \right)^{\!\circ n}\, G \bigg\Vert_{\alpha}\ \ \ge\ \ C_{5}\, | \xi |\, \big\Vert G \big\Vert_{\infty}\; +\; \alpha^{n}\, \big\vert G \big\vert_{\alpha}. \] 
\end{lemma} 
\medskip 

\noindent 
\begin{lemma}[\cite{ps:01}] 
\label{zeropointtwo}
Let $F \in \mathscr{F}_{\alpha} (X, \mathbb{R})$ satisfy the approximability condition, as mentioned in equations \eqref{dioF} and \eqref{dioph} of theorem \eqref{ror1}. Then, there exists positive constants $C_{6},\ C_{7}$ and $C_{8}$ such that 
\[ \bigg\Vert \left( \varrho^{-1} \widetilde{\mathfrak{L}}_{\zeta F} \right)^{\! \circ (2 n R)} \bigg\Vert\ \ \le\ \ C_{6}\, | \xi |\, \left( 1 - \frac{1}{| \xi |^{C_{7}}} \right)^{n - 1},\ \ \ \forall n \ge 1, \] 
where $\left| \xi \right|$ is sufficiently large and $R$ is the greatest integer contained in $C_{8} \log \left| \xi \right|$. 
\end{lemma}   
\medskip 

\noindent 
We conclude this section with a lemma that provides a bound on the operator norm for the iterates of the normalised, but perturbed operator $\widetilde{\mathbb{L}}_{\zeta f}$, for some $f \in \mathscr{F}_{\alpha} (I, \mathbb{R})$. We also include a short proof of the same, for readers' convenience. 
\medskip 

\noindent 
\begin{lemma} 
\label{eightpointthree}
Let $f \in \mathscr{F}_{\alpha}(I, \mathbb{R})$ satisfy the approximability condition, as mentioned in equation \eqref{dioph2} of theorem \eqref{ror2}. Let $\zeta = \kappa + i \xi \in \mathbb{C}$. Then, there exists positive constants $C_{9},\ C_{10}$ and $C_{11}$ such that 
\[ \bigg\Vert \left( e^{-\mathbb{P}(\kappa f)} \widetilde{\mathbb{L}}_{\zeta f} \right)^{\!\circ (2 n R)} \bigg\Vert\ \ \le\ \ C_{9}\ \left\vert \xi \right\vert\ \Bigg( 1 - \frac{1}{\left\vert \xi \right\vert^{C_{10}}} \Bigg)^{n - 1},\ \ \ \forall n \ge 1, \] 
where $\left\vert \xi \right\vert$ is sufficiently large and $R$ is the greatest integer contained in $C_{11} \log \left\vert \xi \right\vert$. 
\end{lemma}

\noindent 
\begin{prooof} 
Given $f \in \mathscr{F}_{\alpha} (I, \mathbb{R})$ that satisfies the approximability condition, it is an easy observation that $Q(f) \in \mathscr{F}_{\alpha} (X, \mathbb{R})$ satisfies the appropriate approximability condition. Thus, lemma \eqref{zeropointtwo} comes to our rescue. Further, from the definitions of the operators as stated in \eqref{ruelleT} and \eqref{ruellesimul}, we obtain that 
\begin{equation} 
\label{landl} 
\bigg\Vert \left( \widetilde{\mathbb{L}}_{\zeta f} \right)^{\!\circ n} \bigg\Vert\ \ \le\ \ \bigg\Vert \left( \widetilde{\mathfrak{L}}_{\zeta Q(f)} \right)^{\!\circ n} \bigg\Vert\ \ \ \ \forall n \ge 1. 
\end{equation} 
\end{prooof} 

\section{Proof of the ergodicity theorems} 
\label{ergproof}

\noindent 
In this section, we prove the ergodicity theorems as stated in \eqref{erg1} for the skew-product setting and \eqref{erg2} for the simultaneous dynamics setting. 
\bigskip 

\noindent 
\begin{proof}[Theorem \eqref{erg1}] 
Let $B \subseteq X$ be a completely $T$-invariant set in the $\sigma$-algebra of $X$ with strictly positive measure, {\it i.e.}, $\mu (B) > 0$. As said earlier in section \eqref{prelims}, a natural candidate for $\mu$ is the product measure of the Bernoulli measure on cylinder sets of $\Sigma_{N}^{+}$ and the Lebesgue measure on open intervals of $I$. Since we have assumed that the set $B$ is completely $T$-invariant, {\it i.e.}, $T^{-1} B\, =\, B$ with strictly positive measure, it is sufficient for us to show that $\mu (B) = 1$, to prove the theorem. 
\medskip 

\noindent 
It is easy to observe that, in its most general form, the set $B$ can be expressed as 
\[ B\ \ =\ \ \bigcup_{j\, \ge\, 1} \big( U_{j} \times V_{j} \big), \]
where $U_{j}$'s are cylinder sets in $\Sigma_{N}^{+}$ and $V_{j}$'s are open subsets of $I$. Let $U \times V \subseteq B$ be a product set in this collection where $U = \big[ v_{1}\, v_{2}\, \cdots\, v_{n} \big]$ and $V \subseteq I$ is an open set. Since $B$ is completely $T$-invariant, we have $U \times V \subseteq T^{-1} B$. Thus, there exists $B'$ in the $\sigma$-algebra of $X$ such that $B' \subseteq B$ and $U \times V \subseteq T^{-1} (B')$. For the smallest such subset $B'$, the countably many possibilities for the form of $B'$ are given by 
\begin{enumerate} 
\item[(0)] $\big[ v_{2}\, v_{3}\, \cdots\, v_{n} \big] \times V'$; 
\item[(1)] $\bigcup\limits_{d\, =\, 1}^{N} \Big( \big[ v_{2}\, v_{3}\, \cdots\, v_{n}\, d \big] \times V_{d}' \Big)$; 
\item[(2)] $\bigcup\limits_{d_{1}\, =\, 1}^{N} \bigcup\limits_{d_{2}\, =\, 1}^{N} \Big( \big[ v_{2}\, v_{3}\, \cdots\, v_{n}\, d_{1}\, d_{2} \big] \times V_{(d_{1}\, d_{2})}' \Big)$; 
\item[($\cdot$)]\ \ $\cdots$. 
\end{enumerate} 
\medskip 

\noindent 
It is clear that in the above enumeration of possibilities, the cylinder sets in all cases starting from (1) onwards are subsumed by the cylinder set in case (0). Suppose we also prove that the appropriate open subsets of $I$ in each of the possibilities starting from case (1) are subsumed by the open subset of $I$ in case (0), it is sufficient to merely work with case (0). For the same purpose, we consider the general case (m), as listed in the above possibilities, namely, 
\begin{enumerate} 
\item[(m)] $\bigcup\limits_{d_{1}\, =\, 1}^{N} \bigcup\limits_{d_{2}\, =\, 1}^{N} \cdots \bigcup\limits_{d_{m}\, =\, 1}^{N}\Big( \big[ v_{2}\, v_{3}\, \cdots\, v_{n}\, d_{1}\, d_{2}\, \cdots\, d_{m} \big] \times V_{(d_{1}\, d_{2}\, \cdots\, d_{m})}' \Big)$. 
\end{enumerate}  
\medskip 

\noindent 
Observe that 
\begin{eqnarray} 
\label{UtimesV} 
U \times V & \subseteq & T^{-1} B' \nonumber \\ 
& = & T^{-1} \Bigg( \bigcup\limits_{d_{1}\, =\, 1}^{N} \bigcup\limits_{d_{2}\, =\, 1}^{N} \cdots \bigcup\limits_{d_{m}\, =\, 1}^{N}\Big( \big[ v_{2}\, v_{3}\, \cdots\, v_{n}\, d_{1}\, d_{2}\, \cdots\, d_{m} \big] \times V_{(d_{1}\, d_{2}\, \cdots\, d_{m})}' \Big) \Bigg) \nonumber \\
& = & \bigcup\limits_{d_{1}\, =\, 1}^{N} \bigcup\limits_{d_{2}\, =\, 1}^{N} \cdots \bigcup\limits_{d_{m}\, =\, 1}^{N}\Bigg( T^{-1} \Big( \big[ v_{2}\, v_{3}\, \cdots\, v_{n}\, d_{1}\, d_{2}\, \cdots\, d_{m} \big] \times V_{(d_{1}\, d_{2}\, \cdots\, d_{m})}' \Big) \Bigg) \nonumber \\ 
& = & \bigcup\limits_{d_{1}\, =\, 1}^{N} \bigcup\limits_{d_{2}\, =\, 1}^{N} \cdots \bigcup\limits_{d_{m}\, =\, 1}^{N} \Bigg( \bigcup_{d\, =\, 1}^{N} \Big( \big[ d\, v_{2}\, v_{3}\, \cdots\, v_{n}\, d_{1}\, d_{2}\, \cdots\, d_{m} \big] \times T_{d}^{-1} V_{(d_{1}\, d_{2}\, \cdots\, d_{m})}' \Big) \Bigg). \nonumber \\ 
& & 
\end{eqnarray} 

\noindent 
However, since $U = \big[ v_{1}\, v_{2}\, \cdots\, v_{n} \big]$, the only possibility in the right hand side of the countable union in equation \eqref{UtimesV} where $U \times V$ can be a subset reduces to 
\[ U \times V\ \ \subseteq\ \ \bigcup\limits_{d_{1}\, =\, 1}^{N} \bigcup\limits_{d_{2}\, =\, 1}^{N} \cdots \bigcup\limits_{d_{m}\, =\, 1}^{N} \Big( \big[ v_{1}\, v_{2}\, v_{3}\, \cdots\, v_{n}\, d_{1}\, d_{2}\, \cdots\, d_{m} \big] \times T_{v_{1}}^{-1} V_{(d_{1}\, d_{2}\, \cdots\, d_{m})}' \Big). \]

\noindent 
The above union is over sets that are disjoint with respect to their first component and thus $V \subseteq T_{v_{1}}^{-1} V_{(d_{1}\, d_{2}\, \cdots\, d_{m})}'$, for any $m$-lettered word $(d_{1}\, d_{2}\, \cdots\, d_{m})$. Now using minimality of the set $B'$, we observe that the set $V_{(d_{1}\, d_{2}\, \cdots\, d_{m})}'$ remains independent of the $m$-lettered word $(d_{1}\, d_{2}\, \cdots\, d_{m})$. Hence, $B' = \big[ v_{2}\, v_{3}\, \cdots\, v_{n} \big] \times V'$, as mentioned in case (0). 
\medskip 

\noindent 
Based on our observation that for every $U \times V = U_{0} \times V_{0} \subseteq B$, there exists a $B_{1} = U_{1} \times V_{1} \subseteq B$ that satisfies $U_{0} \times V_{0} \subseteq T^{-1} \left( U_{1} \times V_{1} \right)$, we obtain a sequence of sets $\left( U_{k} \times V_{k} \right)$ that satisfy, $\left( U_{k - 1} \times V_{k - 1} \right) \subseteq T^{-1} \left( U_{k} \times V_{k} \right)$. However, this process must end in a finite number of steps, precisely $n$, since $U_{0}$ is a cylinder set that fixes only $n$ many positions. Thus, after those finitely many steps, we obtain $\Sigma_{N}^{+} \times V_{n} \subseteq B$ for some open subset $V_{n} \subseteq I$ with strictly positive Lebesgue measure, {\it i.e.}, $\lambda (V_{n}) > 0$. 
\medskip 

\noindent 
Choose the maximal subset $V \subseteq I$ such that $\Sigma_{N}^{+} \times V \subseteq B$, {\it i.e.}, suppose there exists a subset $V' \subseteq I$ such that $\Sigma_{N}^{+} \times V' \subseteq B$ then, $V' \subseteq V$. For such a maximal subset $V \subseteq I$, consider $\big[ w \big] \times V$ for some choice of $w \in \big\{ 1,\, 2,\, \cdots,\, N \big\}$. It is obvious that $\big[ w \big] \times V \subseteq \Sigma_{N}^{+} \times V$. Then, there exists a $V' \subseteq I$ such that 
\begin{equation} 
\label{wtimesV} 
\big[ w \big] \times V\ \ \subseteq\ \ T^{-1} \big( \Sigma_{N}^{+} \times V' \big)\ \ =\ \ \bigcup_{d\, =\, 1}^{N} \Big( \big[ d \big] \times T_{d}^{-1} V' \Big). 
\end{equation} 
Arguing as earlier, we reduce the countable union in the right hand side of equation \eqref{wtimesV} to 
\[ \big[ w \big] \times V\ \ \subseteq\ \ \big[ w \big] \times T_{w}^{-1} V'. \] 
However, owing to the maximality of $V$, we have $V' \subseteq V$ that implies $T_{w}^{-1} V' \subseteq T_{w}^{-1} V$. Thus, 
\[ \big[ w \big] \times V\ \ \subseteq\ \ \big[ w \big] \times T_{w}^{-1} V'\ \ \subseteq \big[ w \big] \times T_{w}^{-1} V. \] 
This implies $V \subseteq T_{w}^{-1} V$ for any choice of $w \in \big\{ 1,\, 2,\, \cdots,\, N \big\}$. However, each of these interval maps preserves the Lebesgue measure, {\it i.e.}, $\lambda (V) = \lambda (T_{w}^{-1} V)$ for all $w \in \big\{ 1,\, 2,\, \cdots,\, N \big\}$. Thus, by eliminating an appropriate set of Lebesgue measure zero from $V$, we obtain $\widetilde{V}$ that satisfies $\widetilde{V} = T_{w}^{-1} \widetilde{V}$ for all $w \in \big\{ 1,\, 2,\, \cdots,\, N \big\}$. Finally, we appeal to the ergodicity of each of these interval maps $T_{w}$ to conclude that the completely $T_{w}$-invariant set $\widetilde{V}$ for all $w$ of strictly positive measure must be of Lebesgue measure $1$. 
\medskip 

\noindent 
Thus, the completely $T$-invariant set $B$ that satisfies $\mu (B) > 0$ has measure $1$, thereby completing the proof of theorem \eqref{erg1}. 
\end{proof} 
\bigskip 

\noindent 
\begin{remark} 
\label{erg1rem} 
The above proof, in fact provides a stronger result than is stated in theorem \eqref{erg1}, namely, for any set $B$ in the $\sigma$-algebra of $X$ that is completely $T$-invariant, meaning $T^{-1} B = B$, we have either $\mu (B) = 0$ or $B$ can be expressed as $\Sigma_{N}^{+} \times \widetilde{V}$ where $\widetilde{V}$ is a set of full Lebesgue measure in $I$, {\it i.e.}, $\lambda \left( \widetilde{V} \right) = 1$. 
\end{remark} 
\medskip 

\noindent 
We now prove the ergodic theorem for simultaneous action of finitely many interval maps, as stated in theorem \eqref{erg2}. 
\bigskip 

\noindent 
\begin{proof}[Theorem \eqref{erg2}] 
For a Lebesgue integrable real-valued function $f$ defined on $I$, define $F := Q(f) \in L^{1} (\mu)$, where $\mu$ is the product measure of the Bernoulli measure on the cylinder sets of $\Sigma_{N}^{+}$ and the Lebesgue measure on the open intervals of $I$. Define 
\[ E\ :=\ \Bigg\{ (w,\, x) \in X\; :\; \liminf_{n\, \rightarrow\, \infty} \frac{1}{n} \sum\limits_{j\, =\, 0}^{n - 1} \mathscr{Q}^{j} (F) (w,\, x)\ =\ \limsup_{n\, \rightarrow\, \infty} \frac{1}{n} \sum\limits_{j\, =\, 0}^{n - 1} \mathscr{Q}^{j} (F) (w,\, x)\ =\ \int\! F\, d \mu \Bigg\}. \] 
One can easily observe that $E$ is a $T$-invariant subset of $X$. Thus, $\mu (E)$ is either zero or one, by theorem \eqref{erg1}. The set of points collected in $E$ is the set of all points in $X$ that satisfies the Birkhoff's pointwise ergodic theorem to the dynamical system $T$ acting on $X$. Hence, $\mu (E) = 1$. Further, from remark \eqref{erg1rem}, we have that $E = \Sigma_{N}^{+} \times E'$, where $\lambda(E') = 1$. Hence, we now have for $\lambda$-almost every $x \in I$ and for every $w \in \Sigma_{N}^{+}$,
\[ \lim_{n\, \rightarrow\, \infty} \frac{1}{n} \sum\limits_{j\, =\, 0}^{n - 1} \mathscr{Q}^{j} (F) (w,\, x)\ \ =\ \ \int\! F\, d \mu. \] 
\medskip 

\noindent 
\begin{claim} 
\label{thereshold}
For a fixed $x_{0} \in E'$, given $\epsilon > 0$, there exists $M_{\epsilon} \in \mathbb{N}$, independent of $w$, such that for all $n \geq M_{\epsilon}$, we have 
\[ \left| \lim_{n\, \rightarrow\, \infty} \frac{1}{n} \sum\limits_{j\, =\, 0}^{n - 1} \mathscr{Q}^{j} (F) (w,\, x_{0})\ -\ \int\! F\, d \mu \right|\ \ <\ \ \epsilon,\ \ \ \forall w \in \Sigma_{N}^{+}. \] 
\end{claim}
\medskip 

\noindent 
We initially prove the theorem, assuming claim \eqref{thereshold} to be true. We shall prove the claim immediately thereafter. From the definition of the composition operators, we have for any arbitrary $\epsilon > 0$, there exists $M_{\epsilon} \in \mathbb{N}$ such that for every $n \geq M_{\epsilon}$, 
\[ \left| \frac{1}{n} f_{w}^{n} (x)\ -\ \int_{0}^{1}\! f\, d \lambda \right|\ \ \le\ \ \epsilon,\ \ \text{for}\ \lambda\text{-a.e.}\ x \in I\ \text{and}\ \forall w = (w_{1}\, w_{2}\, \cdots\, w_{n}). \] 
This implies 
\[ \left| \sum\limits_{w\ :\ |w|\, =\, n} \left( \frac{1}{n} f_{w}^{n} (x)\ -\ \int_{0}^{1}\! f\, d \lambda \right) \right|\ \ \le\ \ N^{n} \epsilon. \] 
Thus, 
\[ \left| \frac{1}{n} \frac{1}{N^{n}} \sum\limits_{w\ :\ |w|\, =\, n} f_{w}^{n} (x)\ -\ \int_{0}^{1}\! f\, d \lambda \right|\ \ \le\ \ \epsilon,\ \ \text{for}\ \ \lambda\text{-a.e.}\ x \in I, \] 
proving the theorem. 
\end{proof} 
\bigskip 

\noindent 
We now complete this section, by proving the claim in \eqref{thereshold}. 
\medskip 

\noindent 
\begin{proof}[Claim \eqref{thereshold}] 
For a fixed $x_{0} \in E'$, consider the sequence $\displaystyle{\left\{ \frac{1}{n} \sum\limits_{j\, =\, 0}^{n - 1} \mathscr{Q}^{j} (F) (w,\, x_{0}) \right\}_{n\, \ge\, 1}}$ of functions defined on the compact space $\Sigma_{N}^{+}$, converging to the constant $\int\! F\, d \mu$. This convergence is uniform over $w \in \Sigma_{N}^{+}$, {\it i.e.}, given any $\epsilon > 0$, there exists $M_{1} = M_{1} (\epsilon)$ such that for all $n > M_{1}$, we have 
\begin{equation} 
\label{kirone}
\left| \frac{1}{n} \sum_{j\, =\, 0}^{n - 1} \mathscr{Q}^{j} (F) (w,\, x_{0})\ -\ \frac{1}{n} \sum_{j\, =\, 0}^{n - 1} \mathscr{Q}^{j} (F) (v,\, x_{0}) \right|\ \ \le\ \ \frac{\epsilon}{2},\ \ \forall w, v \in \Sigma_{N}^{+}. 
\end{equation}
We already know that for some fixed $v \in \Sigma_{N}^{+}$, there exists $M_{2} = M_{2}(\epsilon, v)$ such that for all $n > M_{2}$, we have 
\begin{equation} 
\label{kirtwo}
\left| \frac{1}{n} \sum\limits_{j\, =\, 0}^{n - 1} \mathscr{Q}^{j} (F) (v,\, x_{0})\ -\ \int\! F\, d \mu \right|\ \ \le\ \ \frac{\epsilon}{2}. 
\end{equation}
Taking $M_{3} := \max\{ M_{1},\, M_{2} \}$ where $M_{1}$ and $M_{2}$ are the quantities prescribed by equations \eqref{kirone} and \eqref{kirtwo}, we get for all $n > M_{3}$ and $w \in \Sigma_{N}^{+}$,
\begin{eqnarray*} 
& & \left| \frac{1}{n} \sum_{j\, =\, 0}^{n - 1} \mathscr{Q}^{j} (F) (w,\, x_{0})\ -\ \int\! F\, d \mu \right| \\ 
& \le & \left| \frac{1}{n} \sum_{j\, =\, 0}^{n - 1} \mathscr{Q}^{j} (F) (w,\, x_{0})\ -\ \frac{1}{n} \sum_{j\, =\, 0}^{n - 1} \mathscr{Q}^{j} (F) (v,\, x_{0}) \right|\ +\ \left| \frac{1}{n} \sum_{j\, =\, 0}^{n - 1} \mathscr{Q}^{j} (F) (v,\, x_{0})\ -\ \int\! F\, d \mu \right| \\ 
& \le & \epsilon.
\end{eqnarray*} 
It is clear from the definition that $M_{3}$ is independent of the word $w$ and only dependent on $\epsilon$ and $x_{0}$. 
\end{proof}
\bigskip 

\section{Rates of recurrence} 
\label{rorproof}

\noindent 
In this section, we write the proofs of theorems \eqref{ror1} and \eqref{ror2}. Throughout this section, we work with a fixed $F \in \mathscr{F}_{\alpha} (X, \mathbb{R})$, along with the corresponding normalised operator $\varrho^{-1} \widetilde{\mathfrak{L}}_{F}$ with equilibrium measure $\mu_{F}$. We know that the iterates of this operator obey the result given in lemma \eqref{zeropointtwo}. Further, for some $\zeta = \kappa + i \xi \in \mathbb{C}$, we note that the pressure function associated with the operator $\mathfrak{L}_{\zeta F}$ and its normalised version $\widetilde{\mathfrak{L}}_{\zeta F}$ are one and the same, owing to the definition of pressure and the method of normalisation. 
\medskip 

\noindent 
The following lemma gives an approximation for the eigenvalue of the normalised yet perturbed operator, $\widetilde{\mathfrak{L}}_{\zeta F}$. We urge the reader to observe that the statement of the lemma and hence, its proof, are merely mentioned for a change of variables along the imaginary variable, even though more is true, as one may obtain from, say \cite{ps:94}. 
\medskip 

\noindent 
\begin{lemma} 
\label{morse}
For $\zeta = \kappa + i \xi$, there exists a change of variables $\Upsilon = \Upsilon(\xi)$ such that for $| \xi | < \delta$, we can expand 
\[ e^{\mathfrak{P} \left( \zeta F \right)}\ \ =\ \ e^{\mathfrak{P} \left( \kappa F \right)}\; \Big( 1\; -\; \Upsilon^{2}\; +\; i\, \Theta \left( \Upsilon \right) \Big), \] 
where $\Theta$ is a real-valued function that satisfies $\Theta (\Upsilon) = O ( | \Upsilon |^{3} )$. 
\end{lemma}
\medskip 

\noindent 
\begin{prooof}
By perturbation theory, we know that there exists $\delta > 0$ such that, for $G \in \mathscr{F}_{\alpha}(X, \mathbb{C})$ satisfying $\| G - \kappa F \|_{\alpha} < \delta$, 
\begin{equation} 
\label{analytic}
G\ \ \longmapsto\ \ \mathfrak{P}(G)\ \ \ \ \text{is an analytic map}. 
\end{equation}
 
\noindent 
The analyticity of the above map from the Banach space $\mathscr{F}_{\alpha} (X, \mathbb{C})$ to $\mathbb{C}$, in a neighbourhood of $\kappa F$ assures the existence of a linear map, say $\mathfrak{D} : \mathscr{F}_{\alpha} (X, \mathbb{C}) \longrightarrow \mathbb{C}$ such that 
\[ \lim_{G\, \rightarrow\, \kappa F} \frac{\mathfrak{P} \big( G \big)\ -\ \mathfrak{P} \big( \kappa F \big)\ -\ \mathfrak{D} \big( G - \kappa F \big)}{\| G\ -\ \kappa F \|_{\alpha}}\ \ =\ \ 0, \] 
where $\mathfrak{D}$ is the differential of $\mathfrak{P}$ at $\kappa F$. From equation \eqref{firstderivative}, we have for $\xi \in \mathbb{R}$,
\[ \lim_{\xi \rightarrow 0} \frac{\mathfrak{P} \big( ( \kappa + \xi) F \big)\ -\ \mathfrak{P} \big( \kappa F \big)\ -\ \mathfrak{D} \big( \xi F \big)}{\| \xi F \|_{\alpha}}\ \ =\ \ 0, \] 
for the choice 
\[ \frac{\mathfrak{D} \big( F \big)}{\| F \|_{\alpha}}\ \ =\ \ \left. \frac{d}{d\xi} \mathfrak{P} \big( ( \kappa + \xi ) F \big) \right|_{\xi\, =\, 0}\ \ =\ \ \int\! F\, d \mu_{\kappa F}. \] 
Since $\mathfrak{D}$ is linear, we get
\[ \left. \frac{d}{d \xi} \mathfrak{P} \big( ( \kappa + i \xi ) F \big) \right|_{\xi\, =\, 0}\ \ =\ \ \lim_{\xi\, \rightarrow\, 0} \frac{\mathfrak{D} \big( i \xi F \big)}{\| i \xi F \|_{\alpha}}\ \ =\ \ i \lim_{\xi\, \rightarrow\, 0} \frac{\mathfrak{D} \big( \xi F \big)}{\| \xi F \|_{\alpha}}\ \ =\ \ \int\! i F\, d \mu_{\kappa F}\ \ =\ \ 0. \] 
\medskip 

\noindent 
Similarly from equation \eqref{secondderivative} we have for $\xi \in \mathbb{R}$,
\[ \left. \frac{d^{2}}{d \xi^{2}} \mathfrak{P} \big( ( \kappa + \xi ) F \big) \right|_{\xi\, =\, 0}\ \ =\ \ \lim_{n\, \rightarrow\, \infty} \frac{1}{n} \int\! \Big( F^{n} ((w,\, x)) \Big)^{2}\, d \mu_{\kappa F}. \] 
This implies 
\[ \left. \frac{d^{2}}{d \xi^{2}} \mathfrak{P} \big( ( \kappa + i \xi ) F \big) \right|_{\xi\, =\, 0}\ \ =\ \ \lim_{n\, \rightarrow\, \infty} \frac{1}{n} \int\! \Big( i F^{n} \Big)^{2}\, d \mu_{\kappa F}\ \ =\ \ \lim_{n\, \rightarrow\, \infty} \frac{-1}{n} \int\! \Big( F^{n} \Big)^{2}\, d \mu_{\kappa F}\ \ <\ \ 0. \] 
\medskip 

\noindent 
Since $G \longmapsto \mathfrak{P} (G)$ is an analytic map, as mentioned in \eqref{analytic}, we have the map $\zeta \longmapsto \mathfrak{P} ( \zeta F )$ to be analytic too in a neighbourhood of $\kappa$, where $\zeta = \kappa + i \xi \in \mathbb{C}$. This implies that ${\rm Im} \big( \mathfrak{P} ( \zeta F ) \big)$ is a harmonic function around $\kappa$, {\it i.e.}, for $\xi \in \mathbb{R}$,
\[ \left. \frac{d^{2}}{d \xi^{2}} {\rm Im} \big( \mathfrak{P} ( ( \kappa + \xi ) F ) \big) \right|_{\xi\, =\, 0}\ \ =\ \ \left. \frac{d^{2}}{d \xi^{2}} {\rm Im} \big( \mathfrak{P} ( ( \kappa + i \xi ) F ) \big) \right|_{\xi\, =\, 0}. \] 
We know that for $\xi \in \mathbb{R},\ \mathfrak{P} \big( ( \kappa + \xi ) F \big) \in \mathbb{R}$, {\it i.e.}, 
\[ \left. \frac{d^{2}}{d \xi^{2}} {\rm Im} \big( \mathfrak{P} ( ( \kappa + \xi ) F ) \big) \right|_{\xi\, =\, 0}\ \ =\ \ 0\ \ \ \ \Longrightarrow\ \ \ \ \left. \frac{d^{2}}{d \xi^{2}} {\rm Im} \big( \mathfrak{P} ( ( \kappa + i \xi ) F ) \big) \right|_{\xi\, =\, 0}\ \ =\ \ 0. \] 
\medskip 

\noindent 
Thus we have 
\begin{eqnarray*} 
\left. \frac{d}{d \xi} \mathfrak{P} \big( ( \kappa + i \xi ) F \big) \right|_{\xi\, =\, 0} & = & 0; \\ 
\left. \frac{d^{2}}{d \xi^{2}} {\rm Re} \left( \mathfrak{P} \big( ( \kappa + i \xi ) F \big) \right) \right|_{\xi\, =\, 0} & < & 0; \\ 
\left. \frac{d^{2}}{d \xi^{2}} {\rm Im} \left( \mathfrak{P} \big( ( \kappa + i \xi ) F \big) \right) \right|_{\xi\, =\, 0} & = & 0. 
\end{eqnarray*} 

\noindent 
Thus, the map $\xi \longmapsto {\rm Re} \left( \mathfrak{P} \big( ( \kappa + i \xi ) F \big) \right)$ satisfies the hypothesis of Morse lemma for non degenerate critical points. Thus, similar to lemma (4) in \cite{ps:94} that has been derived in an analogous context as lemma (6) in \cite{ss:07}, we obtain a change of variables $\Upsilon = \Upsilon(\xi)$ for $- \delta < \xi < \delta$ such that 
\[ e^{\mathfrak{P} \big( ( \kappa + i \xi ) F \big)}\ \ =\ \ e^{\mathfrak{P} ( \kappa F )} \big( 1\ -\ \Upsilon^{2}\ +\ i \Theta \left( \Upsilon \right) \big), \] 
where $\Theta \left( \Upsilon ( \xi ) \right) = e^{- \mathfrak{P} ( \kappa F )} {\rm Im} \big( e^{\mathfrak{P} \big( ( \kappa + i \xi ) F \big)} \big)$ and $\Theta \left( \Upsilon \right) = O \big( \left| \Upsilon \right|^{3} \big)$, thus proving the lemma.
\end{prooof}
\medskip 

\noindent 
We now start by considering the left hand side of equation \eqref{ror1eq}, as mentioned in theorem \eqref{ror1}, that measures the cardinality of the set 
\[ \Big\{ (w,\, x)\, \in\, {\rm Fix}_{n} (T)\ :\ a\, \leq\, F^{n}((w,\, x))\, \leq\, b \Big\}. \] 
It is clear that it can be expressed in terms of the indicator function, $\chi_{[a,\, b]}$, {\it i.e.}, 
\[ \# \Big\{ (w,\, x)\, \in\, {\rm Fix}_{n} (T)\ :\ a\, \leq\, F^{n}((w,\, x))\, \leq\, b \Big\}\ \ =\ \ \sum_{(w,\, x)\, \in\, {\rm Fix}_{n} (T)} \chi_{[a,\, b]} \Big( F^{n} ((w,\, x)) \Big). \] 
Since we know that any characteristic function can be approximated by a sequence of smooth functions with compact support under the integral norm, we initially prove a slightly modified result, as stated in proposition \eqref{prop 1}, where the indicator function $\chi_{[a,\, b]}$ is replaced by some smooth function with compact support, say $\tau : \mathbb{R} \longrightarrow \mathbb{R} $. 
\medskip 

\noindent 
\begin{proposition} 
\label{prop 1}
Suppose $F \in \mathscr{F}_{\alpha} (X, \mathbb{R})$ satisfies the hypothesis in theorem \eqref{ror1}, {\it i.e.}, $F$ satisfies the approximability condition and there exists a unique real number $\kappa$ such that $\int\! F\, d \mu_{\kappa F} = 0$. Then, there exists a positive constant $C_{12} > 0$ such that 
\[ \digamma_{\!\! \tau} (n)\ \ :=\ \ \sum_{(w,\, x)\; \in\; {\rm Fix}_{n} (T)} \tau \Big( F^{n} ((w,\, x)) \Big)\ \ \sim\ \ C_{12}\; \frac{e^{n\, \mathfrak{P}(\kappa F)}}{\sqrt{n}}\; \int_{\mathbb{R}}\! \tau (t)\, e^{- \kappa t}\; d t. \] 
\end{proposition}
\medskip 

\noindent 
\begin{prooof}
\noindent 
For $y \in \mathbb{R}$ and $\kappa$ as in the hypothesis, define 
\begin{equation} 
\label{tausubkappa} 
\tau_{\kappa} (y)\ \ :=\ \ \tau (y) e^{-\, \kappa y}, 
\end{equation} 
in order that $\digamma_{\!\! \tau} (n)$ can now be expressed as 
\begin{equation} 
\label{taukappa} 
\digamma_{\!\! \tau} (n)\ \ =\ \ \sum_{(w,\, x)\; \in\; {\rm Fix}_{n} (T)} \tau_{\kappa}\; \Big( F^{n} ((w,\, x)) \Big)\; e^{\kappa F^{n} ((w,\, x))}. 
\end{equation} 
Using inverse Fourier transform and Fubini's theorem, we rewrite equation \eqref{taukappa} as 
\begin{eqnarray*}
\digamma_{\!\! \tau} (n) & = & \sum_{(w,\, x)\; \in\; {\rm Fix}_{n} (T)} \int_{\mathbb{R}}\! \widehat{\tau_{\kappa}} (\xi)\; e^{(\kappa\, +\, i \xi)\; F^{n} ((w,\, x))}\; d \xi \\ 
& = & \int_{\mathbb{R}}\! \widehat{\tau_{\kappa}} (\xi) \sum_{(w,\, x)\; \in\; {\rm Fix}_{n} (T)} e^{(\kappa\, +\, i \xi)\; F^{n} ((w,\, x))}\; d \xi. 
\end{eqnarray*}

\noindent 
By definition, $\tau_{\kappa}$ lives inside a compact support. Hence, by an application of the Paley-Wiener theorem, we deduce that it is sufficient to estimate 
 \[ \sum\limits_{(w,\, x)\; \in\; {\rm Fix}_{n} (T)} e^{(\kappa\, +\, i \xi)\; F^{n}((w,\, x))},\ \ \text{for some}\ \xi\; \in\; \mathbb{R}. \] 
The following lemma, from \cite{dr:73} helps us approximate this sum in terms of the iterates of the appropriate normalised operator such that 
\[ \widetilde{\mathfrak{L}}_{\kappa F} \mathbf{1}\ \ =\ \ \varrho\, \mathbf{1}\ \ \ \ \text{where we recall}\ \ \varrho\ \ =\ \ e^{\mathfrak{P}(\kappa F)}. \] 
\medskip 

\noindent 
\begin{lemma} 
\label{ror1lemma1}
Let $\kappa \in \mathbb{R}$ be the unique real number that satisfies $\int\! F\, d \mu_{\kappa F} = 0$. Then, there exists $0 < \eta < \varrho^{-1}$  such that for every point $(w, x)\; \in\; {\rm Fix}_{n} (T)$, we have 
\[ \sum_{(v,\, y)\; \in\; {\rm Fix}_{n} (T)} e^{(\kappa\, +\, i \xi) F^{n}((v,\, y))}\ \ =\ \ \left( \left( \widetilde{\mathfrak{L}}_{(\kappa\, +\, i \xi) F} \right)^{\!\circ n} \mathbf{1} \right) ((w,\, x))\; \Big( 1\; +\; O \left( n \eta^{n}\, \max\big\{ 1,\; |\xi| \big\} \right) \Big). \] 
\end{lemma}
\medskip 
 
\noindent 
Using the above lemma, we write $\digamma_{\!\! \tau} (n)$ as 
\begin{equation} 
\label{split1} 
\digamma_{\!\! \tau} (n)\ \ =\ \ \int_{\mathbb{R}}\! \widehat{\tau_{\kappa}} (\xi) \left( \left( \widetilde{\mathfrak{L}}_{(\kappa\, +\, i \xi) F} \right)^{\!\circ n} \mathbf{1} \right) ((w,\, x))\ \Big( 1\, +\, O \left( n \eta^{n}\, \max\big\{ 1,\; |\xi| \big\} \right) \Big)\; d \xi. 
\end{equation}
\medskip 

\noindent 
The second term in the above equation is dominated by $(\varrho \eta)^{n}$, which converges to zero faster than any polynomial of $n$. In the remainder of the proof, we estimate the first term of equation \eqref{split1}. 
\medskip 

\noindent 
If $\Xi_{\zeta F} : \mathscr{F}_{\alpha}(X, \mathbb{C}) \longrightarrow \mathscr{F}_{\alpha}(X, \mathbb{C})$ is the one-dimensional eigenprojection associated with $\widetilde{\mathfrak{L}}_{\zeta F}$, for $- \delta < \xi < \delta$, we know by perturbation theory that $\Xi_{\zeta F} (\mathbf{1}) = \mathbf{1} + O \big( \left| \Upsilon \right| \big)$. Thus, by perturbation theory and lemma \eqref{morse}, for $- \delta < \xi  < \delta$ and some $0 < \vartheta < 1$, we obtain 
\begin{eqnarray*} 
\left( \widetilde{\mathfrak{L}}_{\zeta F} \right)^{\!\circ n} \mathbf{1} & = & e^{n \mathfrak{P} ( \zeta F )} \Big( 1 + O \left( \big| \Upsilon \big| \right) \Big) + O ( \vartheta^{n} ) \\ 
& = & e^{n \mathfrak{P} \big( \kappa F \big)} \Big( 1\; -\; \Upsilon^{2}\; +\; i\, \Theta \left( \Upsilon \right) \Big)^{n} \Big( 1 + O \left( \big| \Upsilon \big| \right) \Big) + O ( \vartheta^{n} ). 
\end{eqnarray*} 
\medskip 

\noindent 
The above equation facilitates the splitting of the integral in equation \eqref{split1} into two integrals given by, 
\begin{eqnarray} 
\label{split2} 
\int_{\mathbb{R}}\! \widehat{\tau_{\kappa}} (\xi) \left( \left( \widetilde{\mathfrak{L}}_{(\kappa\, +\, i \xi) F} \right)^{\!\circ n} \mathbf{1} \right) ((w,\, x))\; d \xi & = & \int_{| \xi |\, <\, \delta}\! \widehat{\tau_{\kappa}} (\xi) \left( \left( \widetilde{\mathfrak{L}}_{(\kappa\, +\, i \xi) F} \right)^{\!\circ n} \mathbf{1} \right) ((w,\, x))\; d \xi \nonumber \\ 
& & + \int_{| \xi |\, \ge\, \delta}\! \widehat{\tau_{\kappa}} (\xi) \left( \left( \widetilde{\mathfrak{L}}_{(\kappa\, +\, i \xi) F} \right)^{\!\circ n} \mathbf{1} \right) ((w,\, x))\; d \xi. \nonumber \\ 
& & 
\end{eqnarray} 

\noindent 
We first estimate the first integral in equation \eqref{split2}, using the change in variables from lemma \eqref{morse}. 
\noindent 
\begin{eqnarray} 
\label{split3} 
& & \int_{| \xi |\, <\, \delta}\! \widehat{\tau_{\kappa}} (\xi) \left( \left( \widetilde{\mathfrak{L}}_{(\kappa\, +\, i \xi) F} \right)^{\!\circ n} \mathbf{1} \right) ((w,\, x))\, d \xi \nonumber \\ 
& = & \int_{| \xi |\, <\, \delta} e^{n \mathfrak{P} ( \kappa F )} \Big( 1 - \Upsilon^{2} + i \Theta \left( \Upsilon \right) \Big)^{n} \Big( 1 + O \left( \big| \Upsilon \big| \right) \Big) \Big( \widehat{\tau_{\kappa}} \big( \xi ( \Upsilon ) \big) \Big) \frac{d \xi}{d \Upsilon} d \Upsilon + O ( \vartheta^{n} ) \nonumber \\ 
& = & C_{13}\; \widehat{\tau_{\kappa}} (0)\; e^{n \mathfrak{P} ( \kappa F )} \int_{| \xi |\, <\, \delta} \Big( 1 - \Upsilon^{2} + i \Theta \left( \Upsilon \right) \Big)^{n} \Big( 1 + O \left( \big| \Upsilon \big| \right) \Big)\; d \Upsilon\ +\ O(n^{-1})\ +\ O( \vartheta^{n} ), \nonumber \\ 
& & 
\end{eqnarray} 
where $C_{13} > 0$ is a constant dependent on $\tau_{\kappa}$ and the Jacobian of change of variables. We now apply binomial expansion to the expression in the integral of equation \eqref{split3}, thus splitting it into three parts and define them as follows: 
\begin{eqnarray*} 
I_{0} (n) & := & \int_{- \delta}^{\delta}\! \big( 1\ -\ \Upsilon^{2} \big)^{n}\, d \Upsilon; \\ 
\sum_{j\, =\, 1}^{n} I_{j} (n) & := & \left| \sum_{j\, =\, 1}^{n} \binom{n}{j} \int_{- \delta}^{\delta} \big( 1\ -\ \Upsilon^{2} \big)^{n - j}\; \big( i \Theta (\Upsilon) \big)^{j}\; \left( 1\ +\ O \big( |\Upsilon| \big) \right)\, d \Upsilon \right|; \\ 
J(n) & := & \int_{- \delta}^{\delta}\! \big( 1\ -\ \Upsilon^{2} \big)^{n}\; O \big( |\Upsilon| \big)\, d \Upsilon. 
\end{eqnarray*} 
Using techniques from \cite{ps:94}, we estimate the above integrals to get the following inequalities. 
\begin{eqnarray*} 
I_{0} (n) & = & C_{14} \frac{\Gamma \left( n + 1 \right)}{\Gamma \left( n + 1 + \frac{1}{2} \right)} + O ( ( 1 - \delta^{2} )^{n} ) \\ 
\sum_{j\, =\, 1}^{n} I_{j} (n) & \le & C_{15} \frac{\Gamma \left( n + 1 \right)}{\Gamma \left( n + 1 + \frac{1}{2} \right)}\;\frac{1}{\sqrt{n}} \\ 
| J(n) | & \le & C_{16} \frac{\Gamma(n + 1)}{\Gamma(n + 2)}, 
\end{eqnarray*} 
for some positive constants $C_{14},\ C_{15}$ and $C_{16}$. We now use the identity, 
\[ \lim_{n\, \rightarrow\, \infty} \frac{\Gamma(n + \beta)}{\Gamma(n) n^{\beta}}\ \ =\ \ 1, \] 
to conclude that $I_{0} (n)\ \sim\ \frac{1}{\sqrt{n}} C_{14}$ and that the rest of the terms inside the integral in equation \eqref{split3} converge to zero faster than $\frac{1}{\sqrt{n}}$. 
\medskip 

\noindent 
What remains now is the integral in equation \eqref{split2} over $| \xi | \geq \delta$. For that, we make use of the following lemma, which is a standard result in the theory of Fourier transforms. 
\medskip 

\noindent 
\begin{lemma} 
\label{lemmafourier}
Let $\chi : \mathbb{R} \longrightarrow \mathbb{R}$ be a compactly supported $\mathcal{C}^{r}$ function. Then the Fourier transform $\widehat{\chi} (\xi) = O ( | \xi |^{- r} )$ as $\xi \rightarrow \infty$. 
\end{lemma}
\medskip 

\noindent 
Since $\tau_{\kappa}$ is smooth, we may suppose that $\tau_{\kappa}$ is a compactly supported $\mathcal{C}^{r}$ function, for any arbitrary $r \in \mathbb{N}$. Now using lemma \eqref{zeropointtwo} and lemma \eqref{lemmafourier}, we obtain the following expression. 
\[ \int_{| \xi |\, \ge\, \delta} \widehat{\tau_{\kappa}} ( \xi )\; \Big( \left( \widetilde{\mathfrak{L}}_{(\kappa + i \xi ) F} \right)^{\!\circ n} \mathbf{1} \Big)((w, x))\, d \xi\ \ =\ \ O \left( e^{n \mathfrak{P} \big( \kappa F \big)} \int_\delta^{\infty}\! \left( 1 - \frac{1}{\xi^{C_{7}}} \right)^{\frac{n}{2R}} \xi^{- r}\, d \xi \right), \] 
where $R = [C_{8} \log \xi]$. We now prove that the integral in the right hand side of the above quantity approaches zero faster than $\frac{1}{\sqrt{n}}$ as $n$ goes to $\infty$, using techniques from \cite{ps:94}. In order to achieve the same, we further split the integral into two parts, as 
\[ \int_{\delta}^{\infty} \left( 1\ -\ \frac{1}{\xi^{C_{7}}} \right)^{\frac{n}{2R}} \xi^{- r}\, d \xi\ \ =\ \ \int_{\delta}^{n^{\delta'}} \left( 1\ -\ \frac{1}{\xi^{C_{7}}} \right)^{\frac{n}{2R}} \xi^{- r}\, d \xi\ +\ \int_{n^{\delta'}}^{\infty} \left( 1\ -\ \frac{1}{\xi^{C_{7}}} \right)^{\frac{n}{2R}} \xi^{- r}\, d \xi, \] 
where $\delta < \delta' < \frac{1}{C_{7}}$. Thus, by our choice of $\delta'$, we have $C_{7} \delta' < 1$. Hence, we get the following estimates. The convergence rate of the first term to zero is faster than any polynomial while the second part converges to zero at a polynomial rate dependent on $r$.
\begin{eqnarray} 
\label{split5.1} 
\int_{\delta}^{n^{\delta'}} \left( 1\ -\ \frac{1}{\xi^{C_{7}}} \right)^{\frac{n}{2R}} \xi^{- r}\, d \xi & = & O \left( n^{\delta'} \left( 1\ -\ \frac{1}{n^{C_{7} \delta'}} \right)^{\frac{n}{2 \delta' \beta \log n}} \right) \\ 
\label{split5.2} 
\int_{n^{\delta '}}^{\infty} \left( 1\ -\ \frac{1}{\xi^{C_{7}}} \right)^{\frac{n}{2R}} \xi^{- r}\, d \xi & = & O \left( n^{( 1 - r ) \delta'} \right) 
\end{eqnarray} 
Using the bounds in equations \eqref{split5.1} and \eqref{split5.2}, we now estimate the remaining part of equation \eqref{split2} to be the following.
\[ \int_{| \xi |\, \ge\, \delta} \widehat{\tau_{\kappa}} (\xi) \Big( \widetilde{\mathfrak{L}}_{(\kappa + i \xi) F} \Big)^{\! \circ n}\mathbf{1} ((w, x))\, d \xi\ \ =\ \ O \Big( e^{n \mathfrak{P} ( \kappa F )} n^{(1 - r) \delta'} \Big). \] 
Since $\widehat{\tau_{\kappa}}$ is smooth, we prove that the rate of growth of the second term in equation \eqref{split2} is smaller than $\frac{1}{\sqrt{n}} e^{n \mathfrak{P} ( \kappa F )}$, by considering $r > 1 + \frac{1}{2 \delta'}$. Thus, we obtain the following asymptotic relation. 
\[\digamma_{\!\! \tau}(n)\ \ \sim\ \ C_{12} \frac{e^{n\mathfrak{P} (\kappa F)}}{\sqrt{n}} \widehat{\tau_{\kappa}} (0). \] 
\end{prooof}
\medskip 

\noindent 
\begin{proof}[Theorem \eqref{ror1}] 
The proof of theorem \eqref{ror1} now follows from proposition \eqref{prop 1} where we replace the function $F$ by the characteristic function $\chi_{[a,\, b]}$. 
\end{proof} 
\medskip

\noindent 
We now proceed to prove theorem \eqref{ror2}, using techniques similar to those used in the above proof. Hence, we only highlight the important steps to complete this proof. 
\medskip

\noindent 
\begin{proof}[Theorem \eqref{ror2}] 
As earlier, we begin with a generalisation of the left hand side of the assertion of theorem \eqref{ror2} for a compactly supported function and define $\widetilde{\digamma}_{\!\! \tau} (n)$ as 
\[ \widetilde{\digamma}_{\!\! \tau} (n)\ \ :=\ \ \sum\limits_{w\; :\, |w|\, =\, n} \sum\limits_{x\, \in\, {\rm Fix} T_{w}} \tau ( f^{n}_{w} (x) ). \] 
Clearly, $\widetilde{\digamma}_{\chi_{[a,\, b]}} (n)$ coincides with the expression we need to estimate. By replacing $\tau$ by $\tau_{\kappa}$, as defined in equation \eqref{tausubkappa} and applying Fourier transforms we get the following: 
\[ \widetilde{\digamma}_{\tau} (n)\ \ =\ \ \int_{\mathbb{R}}\! \widehat{\tau_{\kappa}} (\xi) \sum\limits_{w\; :\; |w|\, =\, n} \sum\limits_{x\, \in\, {\rm Fix} (T_{w})} e^{((\kappa + i \xi ) f^{n}_{w} (x))}\, d \xi. \] 
Thus, in order to estimate $\widetilde{\digamma}_{\tau} (n)$, we first estimate the sum inside the integral, namely 
\[ \sum\limits_{w\; :\; |w|\, =\, n} \sum\limits_{x\, \in\, {\rm Fix} (T_{w})} e^{((\kappa + i \xi) f^{n}_{w} (x))}\ \ \text{for}\ \ \xi \in \mathbb{R}. \] 
The following lemma gives a relation between the sum to be estimated and the iterates of the normalised Ruelle operator, $\widetilde{\mathbb{L}}_{(\kappa + i \xi) f}$. By observing the relation between $(\widetilde{\mathbb{L}}_{(\kappa + i \xi) f})^{\circ n}$ and $(\widetilde{\mathfrak{L}}_{(\kappa + i \xi) Q(f)})^{\circ n}$, as prescribed in equation \eqref{landl}, we write the following lemma which is nothing but an easy corollary of lemma \eqref{ror1lemma1}. 
\medskip 

\noindent 
\begin{lemma}
Let $\kappa \in \mathbb{R}$ be the unique real number that satisfies $\int\! f\, d \mathfrak{m}_{\kappa f} = 0$. Then, there exists $0 < \eta < e^{- \mathbb{P}(\kappa f)}$ such that for every point $x$ that satisfies $T_{w}^{n} x = x$ for some $n$-lettered word $w$, we have 
\[ \sum_{v\; :\; |v|\, =\, n} \sum\limits_{y\, \in\, {\rm Fix} (T_{v})} e^{(\kappa\, +\, i \xi) f^{n}_{v} (y)}\ \ =\ \ \left( \left( \widetilde{\mathbb{L}}_{(\kappa\, +\, i \xi) f} \right)^{\!\circ n} \mathbf{1} \right) (x)\; \Big( 1\; +\; O \left( n \eta^{n}\, \max\big\{ 1,\; |\xi| \big\} \right) \Big). \] 
\end{lemma}

\noindent Using the above lemma, we rewrite $\widetilde{\digamma}_{\tau}(n)$ as,
\begin{equation} 
\label{integralform}
\widetilde{\digamma}_{\tau}(n) = \int_{\mathbb{R}} \widehat{\tau}_{\kappa}(\xi) \left( \left((\widetilde{\mathbb{L}}_{(\kappa+i\xi)f})^{\circ n} \mathbf{1} \right)(x)\ \left(1 + O(n \eta^{n} \max \{1, |\xi|\})\right) \right) d\xi.
\end{equation} 
The second term in the above integral is dominated by $(\eta e^{\mathbb{P} (\kappa f)})^{n}$ and thus goes to zero, faster than polynomial of $n$, as $n$ tends to infinity. Using  the definition of $\mathbb{P} ( ( \kappa + i \xi ) f )$, theorem \eqref{rotskew} and lemma \eqref{rotsa}, we conclude that the pressure functions $\mathbb{P} ( ( \kappa + i \xi ) f )$ and $\mathfrak{P} ( ( \kappa + i \xi ) Q(f) )$ coincide. Further, since $\int\! Q(f)\, d \mu_{\kappa Q(f)} = \int\! f\, d \mathfrak{m}_{\kappa f} = 0$, we obtain the following lemma as a corollary of lemma \eqref{morse}.

\begin{lemma}
For $\zeta = \kappa + i \xi$, there exists a change of variables $\Upsilon = \Upsilon(\xi)$ such that for $| \xi | < \delta$, we can expand 
\[ e^{\mathbb{P} \left( \zeta f \right)}\ \ =\ \ e^{\mathbb{P} \left( \kappa f \right)}\; \Big( 1\; -\; \Upsilon^{2}\; +\; i\, \Theta \left( \Upsilon \right) \Big), \] 
where $\Theta$ is a real-valued function that satisfies $\Theta (\Upsilon) = O ( | \Upsilon |^{3} )$. 
\end{lemma} 
\medskip 

\noindent 
By perturbation theory, the one dimensional eigenprojection associated with $\widetilde{\mathbb{L}}_{( \kappa + i \xi ) f}$ is of the form 
\[ \Xi_{( \kappa + i \xi ) f} \mathbf{1}\ \ =\ \ \mathbf{1} + O ( | \Upsilon | ). \] 
Thus for $- \delta < \xi < \delta$ and for some $0 < \vartheta < 1$, we have 
\[ \left( \widetilde{\mathbb{L}}_{( \kappa + i \xi ) f} \right)^{\!\circ n} \mathbf{1}\ \ =\ \ e^{n \mathbb{P} \big( ( \kappa + i \xi ) f \big)} \Big( \mathbf{1} + O \left( \big| \Upsilon \big| \right) \Big) + O ( \vartheta^{n} ). \] 
The integral in the equation \eqref{integralform} can now be split accordingly to get
\begin{eqnarray} 
\label{split integral}
\int_{\mathbb{R}}\! \widehat{\tau_{\kappa}} (\xi) \left( \left( \widetilde{\mathbb{L}}_{(\kappa\, +\, i \xi) f} \right)^{\!\circ n} \mathbf{1} \right) (x)\; d \xi & = & \int_{| \xi |\, <\, \delta}\! \widehat{\tau_{\kappa}} (\xi) \left( \left( \widetilde{\mathbb{L}}_{(\kappa\, +\, i \xi) f} \right)^{\!\circ n} \mathbf{1} \right) (x)\; d \xi \nonumber \\ 
& & + \int_{| \xi |\, \ge\, \delta}\! \widehat{\tau_{\kappa}} (\xi) \left( \left( \widetilde{\mathbb{L}}_{(\kappa\, +\, i \xi) f} \right)^{\!\circ n} \mathbf{1} \right) (x)\; d \xi. \nonumber \\ 
& & 
\end{eqnarray} 
We then apply change of variables to the first integral in the above equation to get
\begin{eqnarray*} 
& & \int_{| \xi |\, <\, \delta}\! \widehat{\tau_{\kappa}} (\xi) \left( \left( \widetilde{\mathbb{L}}_{(\kappa\, +\, i \xi) f} \right)^{\!\circ n} \mathbf{1} \right) (x)\, d \xi \\ 
& = & C_{17}\; \widehat{\tau_{\kappa}} (0)\; e^{n \mathfrak{P} ( \kappa f )} \int_{| \xi |\, <\, \delta} \Big( 1 - \Upsilon^{2} + i \Theta \left( \Upsilon \right) \Big)^{n} \Big( 1 + O \left( \big| \Upsilon \big| \right) \Big)\; d \Upsilon\ +\ O(n^{-1})\ +\ O( \vartheta^{n} ), \\ 
\end{eqnarray*} 
where the expression inside the integral is the same as in equation \eqref{split3}. Thus, the techniques used in the proof of theorem \eqref{ror1} can be used again to obtain 
\[ \int_{| \xi |\, <\, \delta}\! \widehat{\tau_{\kappa}} (\xi) \left( \left( \widetilde{\mathbb{L}}_{(\kappa\, +\, i \xi) f} \right)^{\!\circ n} \mathbf{1} \right) (x)\, d \xi\ \sim\ \frac{1}{\sqrt{n}} C_{18} e^{n \mathbb{P} ( \kappa f)} \widehat{\tau_{\kappa}}(0). \]
\medskip 

\noindent 
Now what remains is to prove that the second integral in equation \eqref{split integral} grows at a rate strictly smaller than $\frac{1}{\sqrt{n}} e^{n \mathbb{P} ( \kappa f )}$, as $n \to \infty$, to complete the proof. This can be achieved again using lemma \eqref{eightpointthree}, as earlier. 
\end{proof}
\bigskip 

\section{Decay of correlations}
\label{docsec} 

\noindent 
In this section, we prove theorems pertaining to the decay of correlations in the skew-product setting and in the setting of simultaneous action of finitely many interval maps, namely theorems \eqref{doc1} and \eqref{doc2}. 
\medskip 

\noindent 
Fix $F \in \mathscr{F}_{\alpha} (X, \mathbb{R})$ and consider the normalised Ruelle operator $\widetilde{\mathfrak{L}}_{F}$ along with its corresponding equilibrium measure $\mu_{F}$. Denote by $\mathscr{F}^{F}_{\alpha} (X, \mathbb{R})$, the set of all $\alpha$-H\"{o}lder continuous functions whose integral with respect to $\mu_{F}$ is zero, {\it i.e.}, 
\[ \mathscr{F}^{F}_{\alpha} (X, \mathbb{R})\ \ :=\ \ \Bigg\{ G \in \mathscr{F}_{\alpha} (X, \mathbb{R})\ :\ \int\! G\, d \mu_{F} = 0 \Bigg\}. \] 
It is easily verifiable that $\mathscr{F}^{F}_{\alpha} (X, \mathbb{R})$ is a subspace of $\mathscr{F}_{\alpha} (X, \mathbb{R})$. Further, the space $\mathscr{F}^{F}_{\alpha} (X, \mathbb{R})$ is preserved by the action of the operator $\widetilde{\mathfrak{L}}_{F}$, {\it i.e.}, $\widetilde{\mathfrak{L}}_{F} : \mathscr{F}^{F}_{\alpha} (X, \mathbb{R}) \longrightarrow \mathscr{F}^{F}_{\alpha} (X, \mathbb{R})$. The Ruelle operator theorem states that the action of $\widetilde{\mathfrak{L}}_{F}$ on $\mathscr{F}^{F}_{\alpha} (X, \mathbb{R})$ has a spectral radius strictly smaller than $\varrho = e^{\mathfrak{P} (F)}$. Equivalently one may consider the operator $\varrho^{-1} \widetilde{\mathfrak{L}}_{F}$ on $\mathscr{F}_{\alpha}^{F} (X, \mathbb{R})$ that has a spectral radius, say $\varrho_{F} < 1$. 
\medskip 

\noindent 
We first state and prove a lemma, that will be useful to prove our main results in this section. 
\medskip 

\noindent 
\begin{lemma} 
\label{lem4.1}
For any $\vartheta \in (\varrho_{F},\, 1)$, there exists a positive constant $C_{19} > 0$ such that 
\[ \bigg\Vert \left( \varrho^{-1} \widetilde{\mathfrak{L}}_{F} \right)^{\!\circ n} G \bigg\Vert_{\alpha}\ \ \le\ \ C_{19}\, \vartheta^{n}\, \big\Vert G \big\Vert_{\alpha}\ \ \forall n \ge 0\ \ \text{and}\ \ \forall G \in \mathscr{F}^{F}_{\alpha} (X, \mathbb{R}). \]
\end{lemma}

\noindent 
\begin{prooof}
Fix a number $\vartheta \in (\varrho_{F},\, 1)$. Choose $\epsilon > 0$ such that $\varrho_{F} + \epsilon < \vartheta$. For this $\epsilon$, there exists $M_{\epsilon} \in \mathbb{Z}_{+}$ such that 
\[ \bigg\Vert \left( \varrho^{-1} \widetilde{\mathfrak{L}}_{F} \right)^{\!\circ n} \bigg\Vert^{\frac{1}{n}}\ \ <\ \ \varrho_{F} + \epsilon,\ \ \forall n \ge M_{\epsilon},\ \ \ \text{since}\ \ \varrho_{F}\ \ = \inf\limits_{n\, \ge\, 1} \bigg\Vert \left( \widetilde{\mathfrak{L}}_{F} \right)^{\!\circ n} \bigg\Vert^{\frac{1}{n}}, \] 
where we only consider the action of $\varrho^{-1} \widetilde{\mathfrak{L}}_{F}$ on $\mathscr{F}_{\alpha}^{F} (X, \mathbb{R})$. This implies 
\[ \bigg\Vert \left( \varrho^{-1} \widetilde{\mathfrak{L}}_{F} \right)^{\!\circ n} G \bigg\Vert_{\alpha}\ \ <\ \ \left( \varrho_{F} + \epsilon \right)^{n}\, \big\Vert G \big\Vert_{\alpha}\ \ <\ \ \vartheta^{n}\, \big\Vert G \big\Vert_{\alpha},\ \ \forall n \ge M_{\epsilon},\ \ \text{and}\ \ \forall G \in \mathscr{F}^{F}_{\alpha} (X, \mathbb{R}). \] 
Further, since $\widetilde{\mathfrak{L}}_{F}$ is a bounded operator, there exists a positive constant $D_{0} > 0$ such that 
\[ \bigg\Vert \left( \varrho^{-1} \widetilde{\mathfrak{L}}_{F} \right)^{\!\circ n} G \bigg\Vert_{\alpha}\ \ \le\ \ D_{0}\, \big\Vert G \big\Vert_{\alpha},\ \ \forall n \ge 1. \] 
Now, an application of the Archimedean property of $\mathbb{R}$ results in finitely many finite constants $D_{1},\, D_{2},\, \cdots,\, D_{M_{\epsilon}}$ that satisfy $D_{n} \vartheta^{n} > D_{0}$ for $1 \le n \le M_{\epsilon}$. Hence, 
\[ \bigg\Vert \left( \varrho^{-1} \widetilde{\mathfrak{L}}_{F} \right)^{\!\circ n} G \bigg\Vert_{\alpha}\ \ \le\ \ D_{n}\, \vartheta^{n}\, \big\Vert G \big\Vert_{\alpha},\ \ \text{for}\ 1 \le n \le M_{\epsilon}. \] 
The result now follows by taking $C_{19} = \max \big\{ 1,\, D_{1},\, D_{2},\, \cdots,\, D_{M_{\epsilon}} \big\}$. 
\end{prooof}
\bigskip 

\noindent 
We are now ready to prove theorem \eqref{doc1} that states the decay of correlations result for the skew-product map $T$. 
\medskip 

\noindent 
\begin{proof}[Theorem \eqref{doc1}] 
For any $G \in \mathscr{F}_{\alpha}^{F} (X, \mathbb{R})$, the left hand side of equation \eqref{doce1} in theorem \eqref{doc1} can be written as 
\begin{eqnarray*} 
\int\! \mathscr{Q}^{n} (G) H\, d \mu_{F}\ -\ \int\! G\, d \mu_{F} \int\! H\, d \mu_{F} & = & \int\! \mathscr{Q}^{n} (G) H\, d \mu_{F}\ -\ \int\! \mathscr{Q}^{n} (G)\, d \mu_{F} \int\! H\, d \mu_{F} \\ 
& = & \int\! \mathscr{Q}^{n} (G) \left( H - \int\! H\, d \mu_{F} \right)\, d \mu_{F}. 
\end{eqnarray*} 
Suppose we denote $\displaystyle{\widetilde{H} := \left( H - \int\! H\, d \mu_{F} \right)}$, then it is easy to see that $\widetilde{H} \in \mathscr{F}^{F}_{\alpha} (X, \mathbb{R})$. Further, on the space of $\mu_{F}$-square integrable real-valued functions defined on $X$ denoted by $L^{2} (\mu_{F})$, the operator $\widetilde{\mathfrak{L}}_{F}$ has a natural extension, with its adjoint given by the operator $\mathscr{Q}$, {\it i.e.}, 
\[ \big\langle \mathscr{Q} \Phi,\, \Psi \big\rangle\ \ =\ \ \big\langle \Phi,\, \widetilde{\mathfrak{L}}_{F} \Psi \big\rangle,\ \ \forall \Phi, \Psi \in L^{2} (\mu_{F}). \] 
Hence,
\[ \int\! \mathscr{Q}^{n} (G) \widetilde{H}\, d \mu_{F}\ \ =\ \ \int\! G \left( \varrho^{-1} \widetilde{\mathfrak{L}}_{F} \right)^{\!\circ n} \widetilde{H}\, d \mu_{F}. \] 
Therefore, 
\[ \left\vert \int\! G \left( \varrho^{-1} \widetilde{ \mathfrak{L}}_{F} \right)^{\!\circ n} \widetilde{H}\, d \mu_{F} \right\vert\ \ \le\ \ \int\! \left\vert G \left( \varrho^{-1} \widetilde{\mathfrak{L}}_{F} \right)^{\!\circ n} \widetilde{H} \right\vert\, d \mu_{F}\ \ \le\ \ \left\Vert G \right\Vert_{2}\, \left\Vert \left( \varrho^{-1} \widetilde{\mathfrak{L}}_{F} \right)^{\!\circ n} \widetilde{H} \right\Vert_{2}. \] 
Further, 
\begin{eqnarray*}
\bigg\Vert \left( \varrho^{-1} \widetilde{\mathfrak{L}}_{F} \right)^{\!\circ n} \widetilde{H} \bigg\Vert_{2} & \le & \bigg\Vert \left( \varrho^{-1} \widetilde{\mathfrak{L}}_{F} \right)^{\!\circ n} \widetilde{H} \bigg\Vert_{\alpha} \\ 
& \le & C_{19}\, \vartheta^{n}\, \left\Vert \widetilde{H} \right\Vert_{\alpha}\ \hspace{+6cm} (\text{using lemma \eqref{lem4.1}}) \\ 
& \le & C_{19}\, \vartheta^{n}\, \left( \big\Vert H \big\Vert_{\alpha}\ +\ \left\vert \int H\, d \mu_{F} \right\vert \right)\ \hspace{+2.7cm} (\text{using definition of}\ \widetilde{H}) \\
& \le & 2 C_{19}\, \vartheta^{n} \big\Vert H \big\Vert_{\alpha}. 
\end{eqnarray*} 
Thus, we obtain the result with the constant  $C_{2}\; =\; 2 C_{19}\, \left\Vert G \right\Vert_{2}\, \left\Vert H \right\Vert_{\alpha}$ to complete the proof of theorem \eqref{doc1}. 
\end{proof}
\bigskip 

\noindent 
To prove theorem \eqref{doc2}, we start by considering the functions $f_{d} = - \log | T_{d}' | \in \mathscr{F}_{\alpha} (I, \mathbb{R})$, for $1 \le d \le N$. By the definition of the Ruelle operator $\mathcal{L}_{f_{d}}^{(d)}$, as stated in equation \eqref{ruelled}, we know that 
\[ \Big( \mathcal{L}_{f_{d}}^{(d)} g \Big) (x)\ \ =\ \ \sum_{T_{d} y\, =\, x} \frac{g(y)}{| T_{d}' (y) |}. \] 
Observe that $\mathcal{L}_{f_{d}}^{(d)} = \widetilde{\mathcal{L}}_{f_{d}}^{(d)}$, {\it i.e.}, the operator $\mathcal{L}_{f_{d}}^{(d)}$ has eigenvalue $1$, with corresponding eigenfunction $\mathbf{1}$. Further, it is evident from Boyarsky and G\'ora (\cite{bg:97}, section (4.3)) that the dual operator $\left( \mathcal{L}_{f_{d}}^{(d)} \right)^{\!*}$ fixes the Lebesgue measure $\lambda$ {\it i.e.}, $\left( \mathcal{L}_{f_{d}}^{(d)} \right)^{\!*} \lambda = \lambda$. Moreover, for every $1 \le d \le N$, the operator $\mathscr{O}_{d}$ defined by $\mathscr{O}_{d} g = g \circ T_{d}$ satisfies $\mathcal{L}_{f_{d}}^{(d)} \mathscr{O}_{d} = {\rm id}$, the identity operator in $\mathcal{C} (I, \mathbb{R})$. 
\medskip 

\noindent 
Denoting by $\mathscr{F}_{\alpha}^{\lambda} (I, \mathbb{R})$, the set of all real-valued $\alpha$-H\"older continuous functions on $I$ whose Lebesgue integral is equal to $0$, {\it i.e.}, 
\begin{equation} 
\label{Falphalambda} 
\mathscr{F}_{\alpha}^{\lambda} (I, \mathbb{R})\ \ :=\ \ \left\{ f \in \mathscr{F}_{\alpha} (I, \mathbb{R})\ :\ \int\! f \, d \lambda = 0 \right\}, 
\end{equation} 
and observing that $\mathcal{L}_{f_{d}}^{(d)}$ preserves the space $\mathscr{F}_{\alpha}^{\lambda} (I, \mathbb{R})$ for all $1 \le d \le N$, we now state a lemma whose proof runs {\it mutatis mutandis} as the proof of lemma \eqref{lem4.1}. We know that the action of $\mathcal{L}_{f_{d}}^{(d)}$ on $\mathscr{F}_{\alpha}^{\lambda} (I, \mathbb{R})$ has a spectral radius, say $\rho_{\lambda}^{(d)} < 1$, owing to theorem \eqref{rotd}. 
\medskip 

\noindent 
\begin{lemma} 
\label{elevenpointtwo}
For any $\vartheta^{(d)} \in (\rho_{\lambda}^{(d)},\, 1)$, there exists a constant $C_{20}^{(d)} > 0$ such that 
\[ \bigg\Vert \left( \mathcal{L}_{f_{d}}^{(d)} \right)^{\!\circ n} g \bigg\Vert_{\alpha}\ \ \le\ \ C_{20}^{(d)}\, \left( \vartheta^{(d)} \right)^{n}\, \big\Vert g \big\Vert_{\alpha}\ \ \forall n \ge 1\ \ \text{and}\ \ \forall g \in \mathscr{F}^{\lambda}_{\alpha} (I, \mathbb{R}),\ \ \text{for}\ \ 1 \le d \le N. \] 
\end{lemma}
\medskip 

\noindent 
We are now thoroughly equipped to prove theorem \eqref{doc2}. 
\medskip 

\noindent 
\begin{proof}[Theorem \eqref{doc2}]
For any $g \in \mathscr{F}_{\alpha}^{\lambda} (I, \mathbb{R})$, the left hand side of equation \eqref{doce2} in theorem \eqref{doc2} can be written as 
\begin{eqnarray*} 
\int\! \mathscr{O}_{w} (g) h\, d \lambda\ -\ \int\! g\, d \lambda \int\! h\, d \lambda & = & \int\! \mathscr{O}_{w} (g) h\, d \lambda\ -\ \int\! \mathscr{O}_{w} (g)\, d \lambda \int\! h\, d \lambda \\ 
& = & \int\! \mathscr{O}_{w} (g) \left( h - \int\! h\, d \lambda \right)\, d \lambda.  
\end{eqnarray*} 
Let $\displaystyle{\widetilde{h} = h - \int\! h\, d \lambda} \in \mathscr{F}_{\alpha}^{\lambda} (I, \mathbb{R})$. Then, since $\mathcal{L}_{f_{d}}^{(d)}$ is the adjoint of $\mathscr{O}_{d}$ in the space of Lebesgue square integrable real-valued functions defined on $I,\ L^2(I,\mathbb{R})$, we have 
\begin{eqnarray*} 
\left\vert \int\! \left( \mathscr{O}_{w} g \right) \widetilde{h}\, d \lambda \right\vert & = & \left\vert \int\! g \left( \mathcal{L}_{f_{w_{n}}}^{(w_{n})} \mathcal{L}_{f_{w_{n - 1}}}^{(w_{n - 1})} \cdots \mathcal{L}_{f_{w_{1}}}^{(w_{1})} \right) \widetilde{h}\, d \lambda \right\vert \\ 
& \le & \left\Vert g \right\Vert_{2} \left\Vert \mathcal{L}_{f_{w_{n}}}^{(w_{n})} \mathcal{L}_{f_{w_{n - 1}}}^{(w_{n - 1})} \cdots \mathcal{L}_{f_{w_{1}}}^{(w_{1})} \widetilde{h} \right\Vert_{2}. 
\end{eqnarray*} 
Now making use of the inequality 
\[ \left\Vert \left( \mathcal{L}_{f_{w_{n}}}^{(w_{n})} \mathcal{L}_{f_{w_{n - 1}}}^{(w_{n - 1})} \cdots \mathcal{L}_{f_{w_{1}}}^{(w_{1})} \right) \widetilde{h} \right\Vert_{2}\ \ \le\ \ \left\Vert \left( \mathcal{L}_{f_{w_{n}}}^{(w_{n})} \mathcal{L}_{f_{w_{n - 1}}}^{(w_{n - 1})} \cdots \mathcal{L}_{f_{w_{1}}}^{(w_{1})} \right) \widetilde{h} \right\Vert_{\alpha}, \] 
and redistributing the operators for $1 \le d \le N$, we obtain 
\begin{eqnarray*} 
\left\Vert \left( \mathcal{L}_{f_{w_{n}}}^{(w_{n})} \mathcal{L}_{f_{w_{n - 1}}}^{(w_{n - 1})} \cdots \mathcal{L}_{f_{w_{1}}}^{(w_{1})} \right) \widetilde{h} \right\Vert_{\alpha} & = & \left\Vert \left( \left( \mathcal{L}_{f_{1}}^{(1)} \right)^{\gamma_{1}} \left( \mathcal{L}_{f_{2}}^{(2)} \right)^{\gamma_{2}} \cdots \left( \mathcal{L}_{f_{N}}^{(N)} \right)^{\gamma_{N}} \right) \widetilde{h} \right\Vert_{\alpha} \\ 
& \le & C_{20}^{(1)}\, C_{20}^{(2)}\, \cdots\, C_{20}^{(N)}\, \left( \vartheta^{(1)} \right)^{\gamma_{1}} \left( \vartheta^{(2)} \right)^{\gamma_{2}} \cdots \left( \vartheta^{(N)} \right)^{\gamma_{N}} \left\Vert \widetilde{h} \right\Vert_{\alpha}, 
\end{eqnarray*} 
appealing to lemma \eqref{elevenpointtwo}. 
\medskip 

\noindent 
Finally, defining $C_{20} := C_{20}^{(1)}\, C_{20}^{(2)}\, \cdots\, C_{20}^{(N)}$ and $\vartheta := \max \left\{ \vartheta^{(1)},\, \vartheta^{(2)},\, \cdots,\, \vartheta^{(N)} \right\}$, we obtain 
\begin{eqnarray*} 
\left\vert \int\! \mathscr{O}_{w} (g) h\, d \lambda\ -\ \int\! g\, d \lambda \int\! h\, d \lambda \right\vert & \le & C_{20} \vartheta^{n} \left\Vert g \right\Vert_{2} \left\Vert \widetilde{h} \right\Vert_{\alpha} \\ 
& \le & 2 C_{20} \vartheta^{n} \left\Vert g \right\Vert_{2} \left\Vert h \right\Vert_{\alpha}, 
\end{eqnarray*} 
thus completing the proof. 
\end{proof}
\bigskip 

\section{Almost sure invariance principles} 
\label{asipsec}

\noindent 
In this section, we prove the almost sure invariance principles as stated in theorems \eqref{asip1} and \eqref{asip2} for both the settings, that we focus in this paper. As in section \eqref{docsec}, we begin by fixing a real-valued H\"{o}lder continuous function $F \in \mathscr{F}_{\alpha} (X, \mathbb{R})$ and considering the corresponding normalised Ruelle operator $\varrho^{-1} \widetilde{\mathfrak{L}}_{F}$ along with its equilibrium measure $\mu_{F}$ and the subspace $\mathscr{F}_{\alpha}^{F} (X, \mathbb{R}) \subseteq \mathscr{F}_{\alpha} (X, \mathbb{R})$. The proof closely follows the method of proof given by Pollicott and Sharp in \cite{ps:02} and Sridharan in \cite{ss:09}.
\medskip 

\noindent 
For any function $G \in \mathscr{F}_{\alpha}^{F} (X, \mathbb{R})$, define 
\[ H\ \ :=\ \ \sum\limits_{n\, \ge \, 1} \left( \varrho^{-1} \widetilde{\mathfrak{L}}_{F} \right)^{\!\circ n} G. \] 
Observe that the infinite series that defines $H$ converges, owing to lemma \eqref{lem4.1}. Then, 
\[  \widetilde{\mathfrak{L}}_{F} \Big( G\, +\, H\, -\, \mathscr{Q} (H) \Big)\ \ =\ \  \widetilde{\mathfrak{L}}_{F} G\, +\,  \widetilde{\mathfrak{L}}_{F} H\, -\, \varrho  H\ \ =\ \ \mathbf{0}. \] 
Thus, defining $\Phi\; :=\; G + H - \mathscr{Q} (H)$, we observe that 
\begin{eqnarray*} 
\Big\vert G^{n} ((w,\, x))\; -\; \Phi^{n} ((w,\, x)) \Big\vert & = & \Big\vert \mathscr{Q}^{n} H ((w,\, x))\; -\; H ((w,\, x)) \Big\vert \\ 
& \le & \Big\vert \mathscr{Q}^{n} H ((w,\, x)) \Big\vert\; +\; \Big\vert H ((w,\, x)) \Big\vert \\ 
& \le & 2 \big\Vert H \big\Vert_{\alpha}. 
\end{eqnarray*} 
\medskip 

\noindent 
Thus, we have proved: 
\begin{lemma}(c.f.\cite{ps:02}, Lemma 2)
For any function $G \in \mathscr{F}_{\alpha}^{F} (X, \mathbb{R})$, there exists a function $H \in \mathscr{F}_{\alpha}^{F} (X, \mathbb{R})$ such that $\Phi = G + \big( H - \mathscr{Q} (H) \big)$ satisfies 
\[ \widetilde{\mathfrak{L}}_{F} \Phi\ \ =\ \ 0\ \ \ \ \text{and}\ \ \ \ G^{n} ((w,\, x))\ \ =\ \ \Phi^{n} ((w,\, x)) + O(1). \] 
\end{lemma}
\medskip 

\noindent 
Given that $\Phi$ and $G$ are cohomologous to each other, we have
\begin{equation} 
\label{eqnvar}
\varsigma(G)^{2}\ \ =\ \ \int\! \Phi ((w,\, x))^{2}\, d \mu_{F}\; +\; 2 \sum_{n\, \ge\, 0} \int\! \Phi ((w,\, x)) \Phi (T^{n} ((w,\, x)))\, d \mu_{F}. 
\end{equation}
Since $\widetilde{\mathfrak{L}}_{f} \Phi = 0$, we obtain 
\begin{eqnarray*} 
\int\! \Phi((w,\, x)) \Phi(T^{n} ((w,\, x)))\, d \mu_{F} & = & \int\! \left( \varrho^{-1} \widetilde{\mathfrak{L}}_{F} \right)^{\!\circ n} \left( \Phi ((w,\, x)) \Phi(T^{n} ((w,\, x))) \right)\, d \mu_{F} \\ 
& = & \int\! \left( \varrho^{-1} \widetilde{\mathfrak{L}}_{F}\right)^{\!\circ n} \Phi((w,\, x)) \Phi((w,\, x))\, d \mu_{F} \\ 
& = & 0. 
\end{eqnarray*} 
Therefore, equation \eqref{eqnvar} becomes 
\[ \varsigma(G)^{2}\ \ =\ \ \int\! \Phi ((w,\, x))^{2}\, d \mu_{F}. \] 

\noindent
Let 
\[ \widehat{X}\ \ :=\ \ \Big\{ (w_{n},\, x_{n})_{n\, \le\, 0} \in X^{- \mathbb{N}}\ :\ T((w_{n - 1},\, x_{n - 1})) = (w_{n},\, x_{n}) \Big\}. \] 
For the purposes of proofs in this section, we fix the following notations. Elements in $\widehat{X}$ will be represented as $\overline{(w,\, x)} = (w_{n},\, x_{n})_{n\, \le\, 0}$. Making use of the canonical projection 
\[ {\rm Pr} : \widehat{X} \longrightarrow X\ \ \ \text{defined by}\ \ \ {\rm Pr} \left( \overline{(w,\, x)} \right)\ =\ (w_{0},\, x_{0}), \] 
we denote and define the natural extension of the map $T : X \longrightarrow X$ on $\widehat{X}$ by 
\[ \widehat{T}\ :\ \widehat{X} \longrightarrow \widehat{X}\ \ \ \text{such that}\ \ \ {\rm Pr} \left( \widehat{T} \left( \overline{(w,\, x)} \right) \right)\ =\ T((w_{0},\, x_{0})). \] 
Given a function $\Phi \in \mathcal{C} (X, \mathbb{R})$, let $\widehat{\Phi}$ be its natural extension on $\widehat{X}$ given by 
\begin{equation} 
\label{widehatPhi} 
\widehat{\Phi} \left( \overline{(w,\, x)} \right)\ \ =\ \ \Phi \left( {\rm Pr} \left (\overline{(w,\, x)} \right) \right)\ \ =\ \ \Phi( (w_{0},\, x_{0}) ). 
\end{equation} 
Thus, the function $F \in \mathscr{F}_{\alpha} (X, \mathbb{R})$ that we fixed in the beginning of this section along with its equilibrium measure $\mu_{F}$ are written as $\widehat{F}$ and $\widehat{\mu_{F}}$ on the space $\widehat{X}$. Since $\mu_{F}$ is a $T$-invariant probability measure on $X$, it is clear that $\widehat{\mu_{F}}$ is a $\widehat{T}$-invariant probability measure on $\widehat{X}$. Suppose $\mathscr{B}$ is a $\sigma$-algebra on $X$, define a sequence of $\sigma$-algebras on $\widehat{X}$ by 
\[ \mathscr{B}_{0}\ \ :=\ \ {\rm Pr}^{- 1} \mathscr{B}\ \ \ \ \text{and}\ \ \ \ \mathscr{B}_{n}\ \ :=\ \ \left( \widehat{T} \right)^{n} \left( \mathscr{B}_{0} \right)\ \ \text{for}\ \ n \in \mathbb{N}. \] 

\noindent
On a probability space $(\Omega, \nu)$, let $\big\{ \mathscr{B}_{n} \big\}_{n\, \ge\, 0}$ be an increasing sequence of $\sigma$-algebras and $\big\{ \Psi_{n} : \Omega \longrightarrow \mathbb{R} \big\}_{n\, \ge\, 0}$ be a collection of functions. Then, $\big\{ \Psi_{n}, \mathscr{B}_{n} \big\}_{n\, \ge\, 0}$ is called an \emph{increasing martingale} if $\Psi_{n}$ is $\mathscr{B}_{n}$-measurable and ${\rm E} \left[ \Psi_{n + 1} \mid \mathscr{B}_{n} \right] = \Psi_{n}$  for $n\geq 0$. 
\medskip

\noindent
Thus, defining 
\begin{eqnarray*} 
\left( \widehat{\Phi} \right)^{n} \left( \overline{(w,\, x)} \right) & := & \widehat{\Phi} \left( \left( \widehat{T} \right)^{- 1} \left( \overline{(w,\, x)} \right) \right)\; +\; \widehat{\Phi} \left( \left( \widehat{T} \right)^{- 2} \left( \overline{(w,\, x)} \right) \right) \\ 
& & \hspace{+4cm} +\; \cdots +\; \widehat{\Phi} \left( \left( \widehat{T} \right)^{- n} \left( \overline{(w,\, x)} \right) \right), 
\end{eqnarray*} 
that captures the $n$-th ergodic sum $\Phi^{n}((w,\, x))$, as defined in equation \eqref{nthergodicsumF}, on the base space, helps us form an increasing martingale on $\widehat{X}$, related to $\Phi^{n}$. 
\medskip 

\noindent 
\begin{lemma}\cite{ps:02}
The sequence $\left\{ \left( \widehat{\Phi} \right)^{n}, \mathscr{B}_{n} \right\}_{n\, \ge\, 1}$ forms an increasing martingale on $\widehat{X}$. 
\end{lemma} 
\medskip 

\noindent 
Before we embark on the proof of theorem \eqref{asip1}, we state the Skorokhod embedding theorem, as in Appendix I of \cite{hh:80}. The statement of this theorem will come in handy, in writing the proof. 
\medskip 

\noindent 
\begin{theorem}[Skorokhod embedding theorem]
Let $\left\{ \widehat{\Psi}_{n}, \mathscr{B}_{n} \right\}_{n\, \ge\, 0}$ be a zero mean and square integrable martingale on $\widehat{X}$. Then, there exists a probability space $( \Omega, \mathcal{A}, \nu)$ that supports a Brownian motion $\mathfrak{B}$ such that $\mathfrak{B}(t)$ has variance $t$, an increasing sequence of $\sigma$-algebras $\big\{ \mathcal{F}_{n} \big\}_{n\, \ge\, 0}$ and a sequence of non negative random variables $\big\{ \mathfrak{X}_{n} \big\}_{n\, \ge\, 1}$ such that if $\mathcal{S}_{0} = 0$ and $\mathcal{S}_{n} = \sum\limits_{j\, =\, 1}^{n} \mathfrak{X}_{j}$ for $n \geq 1$, then 
\begin{enumerate} 
\item $\mathfrak{Y}_{n}\ \ :=\ \ \mathfrak{B} \left( \mathcal{S}_{n} \right)\ \ \stackrel{{\rm d}}{=}\ \ \widehat{\Psi}_{n}$, \\ where $\stackrel{{\rm d}}{=}$ represents equality in distribution, {\it i.e.}, for any Borel measurable set $V$ in $\mathbb{R}$,
\[ \widehat{\mu_{F}} \left( \left\{ \overline{(w,\, x)} \in \widehat{X}\ :\ \widehat{\Psi}_{n} \left( \overline{(w,\, x)} \right) \in V \right\} \right)\ \ =\ \ \nu \left( \big\{ \omega \in \Omega\ :\ \mathfrak{Y}_{n} (\omega) \in V \big\} \right); \] 
\item $\mathfrak{Y}_{n}$ and $\mathcal{S}_{n}$ are $\mathcal{F}_{n}$-measurable; 
\item ${\rm E}\left[ \mathfrak{X}_{n} \mid \mathcal{F}_{n - 1}) \right]\ \ =\ \ {\rm E}\left[ \left( \mathfrak{Y}_{n} - \mathfrak{Y}_{n - 1} \right)^{2} \mid \mathcal{F}_{n - 1} \right],\ \nu$-a.e. for $n \ge 1$. 
\end{enumerate} 
\end{theorem}
\medskip 

\noindent 
We now make use of the Skorokhod embedding theorem and prove theorem \eqref{asip1}. 
\medskip 
 
\noindent 
\begin{proof}[Theorem \eqref{asip1}]
Since $\left( \widehat{\Phi} \right)^{n}$ is a square integrable function with mean zero, we can apply the Skorokhod embedding thoerem. Further, making use of the definition of $\widehat{\Phi}$, as given in equation \eqref{widehatPhi}, we obtain 
\[ \mathfrak{Y}_{n}\ \ \stackrel{{\rm d}}{=}\ \ \left( \widehat{\Phi} \right)^{n}\ \ \stackrel{{\rm d}}{=}\ \ \Phi^{n}. \] 
\medskip 

\noindent 
Thus, in order to complete the proof of theorem \eqref{asip1}, we make the following claim. 
\medskip 

\noindent 
\begin{claim} 
\label{claim2} 
Given any $\delta > 0$, 
\[ \mathfrak{Y}_{n} (\omega)\ \ =\ \ \mathfrak{B} (n) (\omega)\; +\;  O \left( n^{\frac{1}{4}\, +\, \delta} \right)\ \ \forall n \ge 0,\ \ \ \nu\text{-a.e.} \] 
\end{claim} 
\medskip 

\noindent 
Pending proof of claim \eqref{claim2}, it follows from the properties of Brownian motion that
\[ \mathfrak{Y}_{\lfloor t \rfloor}\ \ =\ \ \mathfrak{B} (t)\; +\; O \left( t^{\frac{1}{4}\, +\, \delta} \right)\ \ \forall t \ge 0,\ \ \ \nu\text{-a.e.} \]
This proves the theorem. 
\end{proof}
\medskip 

\noindent 
We now prove our claim \eqref{claim2}. 
\medskip 

\noindent 
\begin{proof}[Claim \eqref{claim2}] 
Since $\mathfrak{Y}_{n} = \mathfrak{B} ( \mathcal{S}_{n} )$, we approximate $\mathcal{S}_{n}$ by $n \varsigma(G)^{2}$, as follows. 
\medskip 

\begin{eqnarray} 
\label{eqn12} 
\mathcal{S}_{n}\; -\; n \varsigma(G)^{2} & = & \sum_{j\, =\, 1}^{n} \Big( \mathfrak{X}_{j}\, -\, {\rm E} \big[ \mathfrak{X}_{j} \mid \mathcal{F}_{j - 1} \big] \Big) \nonumber \\ 
& & +\; \sum_{j\, =1\, }^{n} \Big( {\rm E} \big[ \left( \mathfrak{Y}_{j}\, -\, \mathfrak{Y}_{j - 1} \right)^{2} \mid \mathcal{F}_{j - 1} \big]\; -\; \left( \mathfrak{Y}_{j}\, -\, \mathfrak{Y}_{j - 1} \right)^{2} \Big) \nonumber \\ 
& & +\; \sum_{j\, =\, 1}^{n} \Big( \mathfrak{Y}_{j}\, -\, \mathfrak{Y}_{j - 1} \Big)^{2}\; -\; n \varsigma(G)^{2}. 
\end{eqnarray} 
\medskip 

\noindent 
Given any sequence $\left\{ \widehat{\Psi}_{n} \right\}_{n\, \ge\, 0}$ of functions and an increasing sequence of $\sigma$-algebras $\big\{ \mathcal{F}_{n} \big\}_{n\, \ge\, 0}$ such that $\widehat{\Psi}_{n}$ is $\mathcal{F}_{n}$ measurable for all $n \ge 0$, the sequence defined by 
\[ \left\{ \mathbf{\widehat{\Psi}}_{n}\ \ :=\ \ \sum\limits_{j\, =\, 1}^{n} \left( \widehat{\Psi}_{j}\, -\, {\rm E} \big[ \widehat{\Psi}_{j} \mid \mathcal{F}_{j - 1} \big] \right),\; \mathcal{F}_{n} \right\}_{n\, \ge\, 1} \] 
forms a martingale. Hence, the first and the second terms on the right hand side of equation \eqref{eqn12} are martingales. By the strong law of large numbers for martingales, as can be found in \cite{wf:71}, we can see that for every $\delta > 0$ 
\begin{eqnarray*} 
\sum_{j\, =\, 1}^{n} \Big( \mathfrak{X}_{j}\, -\, {\rm E} \big[ \mathfrak{X}_{j} \mid \mathcal{F}_{j - 1} \big] \Big) & = & O \left( n^{\frac{1}{2}\, +\, \delta} \right); \\ 
\sum_{j\, =1\, }^{n} \Big( {\rm E} \big[ \left( \mathfrak{Y}_{j}\, -\, \mathfrak{Y}_{j - 1} \right)^{2} \mid \mathcal{F}_{j - 1} \big]\; -\; \left( \mathfrak{Y}_{j}\, -\, \mathfrak{Y}_{j - 1} \right)^{2} \Big) & = & O \left( n^{\frac{1}{2}\, +\, \delta} \right). 
\end{eqnarray*} 
We can therefore write equation \eqref{eqn12} as 
\begin{equation} 
\label{eqn13} 
\mathcal{S}_{n}\; -\; n \varsigma(G)^{2}\ \ =\ \ \sum_{j\, =\, 1}^{n} \Big( \mathfrak{Y}_{j}\, -\, \mathfrak{Y}_{j - 1} \Big)^{2}\; -\; n \varsigma(G)^{2}\; +\; O \left( n^{\frac{1}{2}\, +\, \delta} \right). 
\end{equation}
\medskip 

\noindent 
We estimate the sum on the right hand side of equation \eqref{eqn13} with the help of the following series. 
\[ \sum_{j\, \ge\, 1} \frac{1}{j^{\frac{1}{2}\, +\, \delta}} \left[ \Big( \mathfrak{Y}_{j}\, -\, \mathfrak{Y}_{j - 1} \Big)^{2}\; -\; \varsigma(G)^{2} \right]. \] 
\medskip 

\noindent 
\begin{lemma} 
The following integrals are equal. 
\begin{eqnarray*} 
\mathfrak{I} (\delta) & := & \int\! \left( \sum_{j\, \ge\, 1} \frac{1}{j^{\frac{1}{2}\, +\, \delta}} \left[ \Big( \mathfrak{Y}_{j}\, -\, \mathfrak{Y}_{j - 1} \Big)^{2}\; -\; \varsigma(G)^{2} \right] \right)^{2}\, d \nu \\ 
& = & \int\! \left( \sum_{j\, \ge\, 1} \frac{1}{j^{\frac{1}{2}\, +\, \delta}} \left[ \Phi \left( T^{j} ((w, x)) \right)^{2}\, -\, \int\! \Phi^{2}\, d \mu_{F} \right] \right)^{2}\, d \mu_{F}. 
\end{eqnarray*} 
\end{lemma} 
\medskip 

\noindent 
\begin{prooof}
We already have that $\mathfrak{Y}_{n} \stackrel{{\rm d}}{=} \Phi^{n}$. Thus, from a proposition of Brieman, L. as in \cite{lb:68}, we deduce that for any measurable function $\Theta : \mathbb{R}^{\mathbb{N}} \longrightarrow \mathbb{R}$,  
\[ \int\! \Theta \Big( \big( \mathfrak{Y}_{j} (\omega) \big)_{j\, =\, 0}^{\infty} \Big)\, d \nu\ \ =\ \ \int\! \Theta \Big( \big( \Phi^{j} ((w, x)) \big)_{j\, =\, 0}^{\infty} \Big)\, d \mu_{F}. \] 
The result follows from an appropriate choice of the function $\Theta$, say 
\[ \Theta \Big( \big( y_{j} \big)_{j\, \ge\, 0} \Big)\ \ =\ \ \left( \sum_{j\, \ge\, 1} \frac{1}{j^{\frac{1}{2}\, +\, \delta}} \left[ \left[y_{j + 1}\, -\, y_{j} \right]^{2}\; -\; \int\! \Phi^{2}\, d \mu \right] \right)^{2}. \] 
\end{prooof} 
\medskip 

\noindent 
A simple calculation now yields that for any $\delta > 0,\ \mathfrak{I} (\delta) < \infty$. Hence, 
\[ \sum_{j\, \ge\, 1} \frac{1}{j^{\frac{1}{2}\, +\, \delta}} \left[ \Big( \mathfrak{Y}_{j}\, -\, \mathfrak{Y}_{j - 1} \Big)^{2}\; -\; \varsigma(G)^{2} \right]\ \ <\ \ \infty,\ \ \nu\text{-a.e.} \] 
Applying the Kronecker lemma as in \cite{hh:80}, we deduce that 
\begin{equation} 
\label{eqn18} 
\sum_{j\, =\, 1}^{n} \left[ \Big( \mathfrak{Y}_{j}\, -\, \mathfrak{Y}_{j - 1} \Big)^{2}\; -\; n \varsigma(G)^{2} \right]\ \ =\ \ O \left( n^{\frac{1}{2}\, +\, \delta} \right). 
\end{equation}
Thus, from equations \eqref{eqn13}) and \eqref{eqn18}, we have $\mathcal{S}_{n}\; -\; n \varsigma(G)^{2}\ =\ O \left( n^{\frac{1}{2}\, +\, \delta} \right),\ \nu$-a.e. 
\medskip 

\noindent 
Finally, defining $\widetilde{\mathfrak{B}} (t) := \mathfrak{B} ( t \varsigma(G)^{2} )$, we have for $n \ge 0$, 
\[ \mathfrak{B} ( \mathcal{S}_{n} )\ \ =\ \ \mathfrak{B} \left( n \varsigma(G)^{2} \right)\; +\; O \left( n^{\frac{1}{4}\, +\, \delta} \right)\ \ =\ \ \widetilde{\mathfrak{B}} (n)\; +\; O \left( n^{\frac{1}{4}\, +\, \delta} \right),\ \ \nu\text{-a.e.} \] 
This proves the equation in claim \eqref{claim2}, namely, 
\[ \mathfrak{Y}_{n}\ \ =\ \ \widetilde{\mathfrak{B}} (n)\; +\; O \left( n^{\frac{1}{4}\, +\, \delta} \right),\ \ \nu\text{-a.e.} \] 
\end{proof} 
\medskip 
  
\noindent 
We now proceed to prove the next theorem on almost sure invariance principles for simultaneous action of the interval maps, as stated in theorem \eqref{asip2}. We draw motivation from the proof of a similar result in an article by Haydn, Nicol, T\"or\"ok and Vaienti in \cite{hntv:17} and achieve a better bound. We first state a theorem due to Cuny and Merl\'evede as in \cite{cm:15}, that would be helpful in our proof. 
\medskip 

\noindent 
\begin{theorem}\cite{cm:15} 
\label{CM 2.3} 
Let $\big\{ U_{n} \big\}_{n\, \ge\, 0}$ be a sequence of square integrable random variables adapted to some non-increasing sequence of $\sigma$-algebras $\big\{ \mathscr{A}_{n} \big\}_{n\, \ge\, 0}$ on $\mathbb{R}$. Assume that 
\[ {\rm E} \big[ U_{n} \mid \mathscr{A}_{n + 1} \big]\ \ =\ \ 0\ \text{a.s.};\ \ \ \ \varsigma_{n}^{2}\ \ =\ \ \sum_{k\, =\, 0}^{n - 1} {\rm E} \big[ U_{k}^{2} \big]\ \ \to\ \ \infty;\ \ \ \ \sup\limits_{n\, \ge\, 0} {\rm E} \big[ U_{n}^{2} \big]\ \ <\ \ \infty. \] 
Let $\big\{ a_{n} \big\}_{n\, \ge\, 0}$ be a non-decreasing sequence of positive numbers such that 
\[ \left\{ \frac{a_{n}}{\varsigma_{n}} \right\}_{n\, \ge\, 0}\ \ \text{is non-decreasing}\ \ \ \text{and}\ \ \ \left\{ \frac{a_{n}}{\varsigma_{n}^{2}} \right\}_{n\, \ge\, 0}\ \ \text{is non-increasing}. \] 
Further, assume that
\begin{enumerate}
\item $\sum\limits_{k\, =\, 0}^{n - 1} \Big( {\rm E} \big[ U_{k}^{2} \mid \mathscr{A}_{k + 1} \big]\; -\; {\rm E} \big[ U_{k}^{2} \big] \Big)\ \ =\ \ o(a_{n}),\ \lambda$-a.s.; 
\item $\sum\limits_{n\, \ge\, 0} a_{n}^{- r} {\rm E} \big[ |U_{n}|^{2r} \big]\ \ <\ \ \infty$ for some $1 \le r \le 2$. 
\end{enumerate} 
Then enlarging our probability space, if necessary, it is possible to find a sequence $\big\{ \mathcal{U}_{n} \big\}_{n\, \ge\, 0}$ of independent centered Gaussian variables with ${\rm E} \big[ \mathcal{U}_{n}^{2} \big] = {\rm E} \big[ U_{n}^{2} \big]$ such that 
\[ \sup_{0\, \le\, k\, \le\, n - 1} \left| \sum_{j\, =\, 0}^{k} U_{j}\; -\; \sum_{j\, =\, 0}^{k} \mathcal{U}_{j} \right|\ \ =\ \ o \left( \left[ a_{n} \left( \left| \log \left( \frac{\varsigma_{n}^{2}}{a_{n}} \right) \right|\; +\; \log \log a_{n} \right) \right]^{\frac{1}{2}} \right),\ \ \ \lambda\text{-a.s.}  \] 
\end{theorem}
\medskip
 
\noindent 
Note that the assertion of theorem \eqref{CM 2.3} can be rewritten by considering another probability space $(\Omega, \mathscr{A}, \nu)$ and a sequence of random variables, say $\big\{ \mathcal{V}_{n} \big\}_{n\, \ge\, 0}$ such that $U_{n} \stackrel{{\rm d}}{=} \mathcal{V}_{n}$. Then, 
\[ \sup_{0\, \le\, k\, \le\, n - 1} \left| \sum_{j\, =\, 0}^{k} \mathcal{V}_{j}\; -\; \sum_{j\, =\, 0}^{k} \mathcal{U}_{j} \right|\ \ =\ \ o \left( \left[ a_{n} \left( \left| \log \left( \frac{\varsigma_{n}^{2}}{a_{n}} \right) \right|\; +\; \log \log a_{n} \right) \right]^{\frac{1}{2}} \right),\ \ \ \nu\text{-a.s.} \] 
\medskip

\noindent	
We will now prove theorem \eqref{asip2}. 
\medskip 

\noindent 
\begin{proof}[Theorem \eqref{asip2}] 
Recall the definition of the space $\mathscr{F}_{\alpha}^{\lambda} (I, \mathbb{R})$ from equation \eqref{Falphalambda}, 
\[ \mathscr{F}_{\alpha}^{\lambda} (I, \mathbb{R})\ \ :=\ \ \left\{ f \in \mathscr{F}_{\alpha} (I, \mathbb{R})\ :\ \int\! f \, d \lambda = 0 \right\}, \] 
and the property that for $f_{d} = - \log | T_{d}' |$, the operator $\mathcal{L}_{f_{d}}^{(d)}$ preserves the space $\mathscr{F}_{\alpha}^{\lambda} (I, \mathbb{R})$ for all $1 \le d \le N$. 
\medskip 

\noindent 
Let $g \in \mathscr{F}_{\alpha}^{\lambda} (I, \mathbb{R})$ and $w \in \Sigma_{N}^{+}$. Then, define a sequence of $\sigma$-algebras 
\[ \mathscr{B}_{w}^{n}\ \ :=\ \ \left( \mathscr{O}_{(w_{1}\, w_{2}\, \cdots\, w_{n})} \right)^{- 1} \mathscr{B}\ \ \ \ \text{for}\ \ n \ge 0, \] 
where $\mathscr{B}$ is the Borel $\sigma$-algebra on $I$. 
\medskip
	
\noindent 
Suppose for all $n \ge 1$, we denote by $\mathfrak{g}_{w}^{n}$, the sum 
\[ \mathfrak{g}_{w}^{n}\ \ :=\ \ \Big( \mathcal{L}_{f_{w_{n}}}^{(w_{n})} \Big) g\; +\; \Big( \mathcal{L}_{f_{w_{n}}}^{(w_{n})} \mathcal{L}_{f_{w_{n - 1}}}^{(w_{n - 1})} \Big) g\; +\; \cdots\; +\; \Big( \mathcal{L}_{f_{w_{n}}}^{(w_{n})} \mathcal{L}_{f_{w_{n - 1}}}^{(w_{n - 1})} \cdots \mathcal{L}_{f_{w_{1}}}^{(w_{1})} \Big) g \] 
and $\mathfrak{g}_{w}^{0} = 0$ for all $w \in \Sigma_{N}^{+}$. It is easy to see that 
\[ \mathcal{L}_{f_{w_{n + 1}}}^{(w_{n + 1})} \mathbb{g}_{w}^{n}\ \ =\ \ 0,\ \ \ \ \text{where}\ \ \ \ \mathbb{g}_{w}^{n}\ \ =\ \ g\; +\; \mathfrak{g}_{w}^{n}\; -\; T_{w_{n + 1}} \mathfrak{g}_{w}^{n + 1}. \] 
Defining $\mathbb{h}_{w}^{n} = \mathscr{O}_{(w_{1}\, w_{2}\, \cdots\, w_{n})} \left( \mathbb{g}_{w}^{n} \right)$, one can observe that $\mathbb{h}_{w}^{n}$ agrees with the definition of a reverse martingale difference sequence for the sequence of $\sigma$-algebras $\mathscr{B}_{w}^{n}$, as defined in Conze and Raugi \cite{cr:07}, as given below. 
\medskip 

\noindent 
Given a sequence of random variables $\big\{ X_{n} \big\}_{n\, \in\, \mathbb{N}}$ adapted to a non-increasing sequence of $\sigma$- algebras $\big\{ \mathscr{A}_{n} \big\}_{n\, \in\, \mathbb{N}},\ \big\{ X_{n},\, \mathscr{A}_{n} \big\}_{n\, \in\, \mathbb{N}}$ is a \emph{reverse martingale} or equivalently, $\big\{ X_{n} \big\}_{n\, \in\, \mathbb{N}}$ is a reverse martingale adapted to $\big\{ \mathscr{A}_{n} \big\}_{n\, \in\, \mathbb{N}}$ if $\bigg\{ \widetilde{X}_{n},\, \widetilde{\mathscr{A}}_{n} \bigg\}_{n\, \le\, -1}$ forms a martingale, where $\widetilde{X}_{n} = X _{- n}$ and $\widetilde{\mathscr{A}}_{n} = \mathscr{A}_{- n}$ for each $n \in - \mathbb{N}$. 
\medskip 

\noindent 
Now, 
\begin{eqnarray*} 
\sum_{k\, =\, 0}^{n - 1} \mathbb{h}_{w}^{k} & = & \sum_{k\, =\, 0}^{n - 1} \left( \mathscr{O}_{(w_{1}\, w_{2}\, \cdots\, w_{k})} g\, +\, \mathscr{O}_{(w_{1}\, w_{2}\, \cdots\, w_{k})} \mathfrak{g}_{w}^{k}\, -\, \mathscr{O}_{(w_{1}\, w_{2}\, \cdots\, w_{k + 1})} \mathfrak{g}_{w}^{k + 1} \right) \\ 
& = & \sum_{k\, =\, 0}^{n - 1} \mathscr{O}_{(w_{1}\, w_{2}\, \cdots\, w_{k})} g\, -\, \mathscr{O}_{(w_{1}\, w_{2}\, \cdots\, w_{n})} \mathfrak{g}_{w}^{n}. 
\end{eqnarray*} 
Further, $\left\Vert \mathfrak{g}_{w}^{n} \right\Vert_{\alpha}$ is uniformly bounded. Hence, 
\begin{eqnarray*} 
{\rm E} \left[ \left( \sum_{k\, =\, 0}^{n - 1} \mathbb{h}_{w}^{k} \right)^{2} \right] & = & {\rm E} \left[ \left( \sum_{k\, =\, 0}^{n - 1} \mathscr{O}_{(w_{1}\, w_{2}\, \cdots\, w_{k})} g \right)^{2} \right]\; +\; {\rm E} \big[ \left( \mathscr{O}_{(w_{1}\, w_{2}\, \cdots\, w_{n})} \mathfrak{g}_{w}^{n} \right)^{2} \big] \\ 
& & \hspace{+2cm} -\; 2 {\rm E} \left[ \left( \sum_{k\, =\, 0}^{n - 1} \mathscr{O}_{(w_{1}\, w_{2}\, \cdots\, w_{k})} g \right) \big( \mathscr{O}_{(w_{1}\, w_{2}\, \cdots\, w_{n})} \mathfrak{g}_{w}^{n} \big) \right] \\ 
& = & \left( \varsigma_{w}^{(n)} (g) \right)^{2}\; +\; o \left( \left( \varsigma_{w}^{(n)} (g) \right) \right),  
\end{eqnarray*} 
where we recall the definition of $\varsigma_{w}^{(n)} (g)$ from equation \eqref{varsigmawng} as 
\[ \left( \varsigma_{w}^{(n)} (g) \right)^{2}\ \ =\ \ \int\! \left( g_{w}^{n} \right)^{2}\, d \lambda. \] 
\medskip 

\noindent 
Haydn, Nicol, T\"or\"ok and Vaienti in \cite{hntv:17} show us that ${\rm E} \big[ \mathbb{h}_{w}^{j} \mathbb{h}_{w}^{k} \big] = 0$, for $j \ne k$ and therefore 
\[ \sum_{k\, =\, 0}^{n - 1} {\rm E} \big[ \left( \mathbb{h}_{w}^{k} \right)^{2} \big]\ \ =\ \ {\rm E} \left[ \left( \sum_{k\, =\, 0}^{n - 1} \mathbb{h}_{w}^{k} \right)^{2} \right]\ \ =\ \ \left( \varsigma_{w}^{(n)} (g) \right)^{2}\; +\; o \left( \left( \varsigma_{w}^{(n)} (g) \right) \right), \] 
which implies, $\sum\limits_{k\, =\, 0}^{n - 1} {\rm E} \big[\left( \mathbb{h}_{w}^{k} \right)^{2} \big] \to \infty$. 
\medskip 

\noindent 
Thus, we have constructed a sequence of square integrable random variables $\left\{ \mathbb{h}_{w}^{n} \right\}_{n\, \ge\, 0}$ adapted to a non-increasing sequence of $\sigma$-algebras $\big\{ \mathscr{B}_{w}^{n} \big\}_{n\, \ge\, 0}$ that satisfies 
\[ {\rm E} \big[ \mathbb{h}_{w}^{n} \mid \mathscr{B}_{w}^{n + 1} \big]\ \ =\ \ 0\ \text{a.s.};\ \ \ \ \sum_{k\, =\, 0}^{n - 1} {\rm E} \big[ \left( \mathbb{h}_{w}^{k} \right)^{2} \big]\ \ \to\ \ \infty;\ \ \ \ \sup\limits_{n\, \ge\, 0} {\rm E} \big[ \left( \mathbb{h}_{w}^{n} \right)^{2} \big]\ \ <\ \ \infty. \] 
Further, defining a sequence $\left\{ a_{n} := \left( \varsigma_{w}^{(n)} (g) \right)^{1\, +\, \epsilon} \right\}_{n\, \ge\, 0}$ for some sufficiently small $\epsilon > 0$, we observe that the sequences satisfy  
\[ \left\{ \frac{a_{n}}{\left( \varsigma_{w}^{(n)} (g) \right)} \right\}_{n\, \ge\, 0}\ \ \text{is non-decreasing}\ \ \ \text{and}\ \ \ \left\{ \frac{a_{n}}{\left( \varsigma_{w}^{(n)} (g) \right)^{2}} \right\}_{n\, \ge\, 0}\ \ \text{is non-increasing}. \] 
\medskip 

\noindent 
Thus, in order to appeal to theorem \eqref{CM 2.3} and exploit the assertions there, we only need to verify the two enumerated assumptions in the statement. We will, for now take the relevant assumptions to be true and proceed to complete the proof of theorem \eqref{asip2}. Once the proof is complete, we will complete the verifications of the enumerated statements. 
\medskip 

\noindent 
By theorem \eqref{CM 2.3}, we have that there exist sequences $\big\{ \mathcal{Y}_{w}^{n} \big\}_{n\, \ge\, 0}$ and $\big\{ Z_{w}^{n} \big\}_{n\, \ge\, 0}$ such that 
\begin{eqnarray*} 
\sup_{0\, \le\, k\, \le\, n - 1} \left| \sum_{j\, =\, 0}^{k} \mathcal{Y}_{w}^{j}\; -\; \sum_{j\, =\, 0}^{k} Z_{w}^{j} \right| & = & o \left( \left[ a_{n} \left( \left| \log \left( \frac{\varsigma_{n}^{2}}{a_{n}} \right) \right|\; +\; \log \log a_{n} \right) \right]^{\frac{1}{2}} \right) \\ 
& = & o \left( \left( n^{\frac{1}{2}\, +\, \epsilon} \left( \left| \log \left( n^{\frac{1}{2}\, -\, \epsilon} \right) \right|\; +\; \log \log n^{\frac{1}{2}\, +\, \epsilon} \right) \right)^{\frac{1}{2}} \right) \\ 
& = & O \left( n^{\frac{1}{4}\, +\, \delta} \right),\ \ \ \ \nu\text{-a.s., for some}\ \delta > 0.  
\end{eqnarray*} 
\medskip 

\noindent 
Further, by the result due to Cuny and Merlevede \cite{cm:15}, we know that	
\[ \sum_{j\, =\, 0}^{n - 1} {\rm E} \big[ \left( Z_{w}^{j} \right)^{2} \big]\ \ =\ \ \left( \varsigma_{w}^{(n)} (g) \right)^{2}\; +\; O \left( \varsigma_{w}^{(n)} (g) \right)\ \ =\ \ \left( \varsigma_{w}^{(n)} (g) \right)^{2}\; +\; o \left( n^{\frac{1}{2}\; +\; \delta'} \right), \] 
for some $\delta' > 0$. Hence, we can replace the random variables with a standard Brownian motion $\big\{ \mathfrak{B}^{*} (t) \big\}_{t\, \ge\, 0}$ such that 
\[ \sup\limits_{0\, \le\, k\, \le\, n - 1} \left| \sum_{j\, =\, 0}^{k} Z_{w}^{j}\; -\; \mathfrak{B}^{*} \left( \left( \varsigma_{w}^{(k)} (g) \right)^{2}\; +\; o \left( k^{\frac{1}{2}\, +\, \delta'} \right) \right) \right|\ \ =\ \ 0,\ \ \ \nu\text{-a.s.} \] 
which implies that
\[ \sup\limits_{0\, \le\, k\, \le\, n - 1} \left| \sum_{j\, =\, 0}^{k} Z_{w}^{j}\; -\; \mathfrak{B}^{*} \left( \left( \varsigma_{w}^{(k)} (g) \right)^{2} \right) \right|\ \ =\ \ o \left( n^{\frac{1}{4}\, +\, \delta'} \right),\ \ \ \nu\text{-a.s.} \] 
Therefore, if we replace  the independent centered Gaussian variables with the standard Brownian motion, we get 
\[ \sum_{j\, =\, 0}^{n - 1} \mathcal{Y}_{w}^{j}\; -\; \mathfrak{B}^{*} \left( \left( \varsigma_{w}^{(n)} (g) \right)^{2} \right)\ \ =\ \ O \left( n^{\frac{1}{4}\, +\, \delta} \right),\ \ \ \ \nu\text{-a.s.} \] 
Further, it is easy to see that there exist a sequence of random variables $\big\{ Y_{w}^{n} \big\}_{n\, \ge\, 0}$ such that 
\[ \left\vert Y_{w}^{n}\; -\; \sum_{j\, =\, 0}^{n} \mathcal{Y}_{w}^{j} \right\vert\ \ =\ \ O(1), \] 
and $g_{w}^{n}$ and $Y_{w}^{n}$ are equal in distribution, thus proving theorem \eqref{asip2}. 
\end{proof} 
\medskip 

\noindent 
We now complete the verifications of the enumerated conditions in theorem \eqref{CM 2.3}. 
\medskip 

\noindent 
\begin{claim} 
\label{cond1} 
The first of the enumerated condition in theorem \eqref{CM 2.3} looks like 
\[ \sum_{k\, =\, 0}^{n - 1} \left( {\rm E} \left[ \left( \mathbb{h}_{w}^{k} \right)^{2} \mid \mathscr{B}_{w}^{k + 1} \right]\; -\; {\rm E} \left[ \left( \mathbb{h}_{w}^{k} \right)^{2} \right] \right)\ \ =\ \ o(a_{n}). \] 
\end{claim} 
\medskip 

\noindent 
\begin{claim} 
\label{cond2} 
The second of the enumerated conditions in theorem \eqref{CM 2.3} looks like 
\[ \sum\limits_{n\, \ge\, 0} a_{n}^{- r} {\rm E} \big[ \left| \mathbb{h}_{w}^{n} \right|^{2r} \big]\ \ <\ \ \infty,\ \ \text{for some}\ 1 \le r \le 2. \] 
\end{claim} 
\medskip 

\noindent 
\begin{proof}[Claim \eqref{cond1}]
From Conze and Raugi \cite{cr:07}, we get that
\[ {\rm E} \big[ \left( \mathbb{h}_{w}^{n} \right)^{2} \mid \mathscr{B}_{w}^{n + 1} \big]\ \ =\ \ \mathscr{O}_{(w_{1}\, w_{2}\, \cdots\, w_{n})} \left( \mathcal{L}_{f_{w_{n + 1}}}^{(w_{n + 1})} \left( \mathbb{g}_{w}^{n} \right)^{2} \right)\ \ \text{and} \]
\[ \int\! \left| \sum_{k\, =\, 0}^{n - 1} \left( {\rm E} \left[ \left( \mathbb{h}_{w}^{k} \right)^{2} \mid \mathscr{B}_{w}^{k + 1} \right]\; -\; {\rm E} \left[ \left( \mathbb{h}_{w}^{k} \right)^{2} \right] \right) \right|^{2}\, d \lambda\ \ \le\ \ C_{21} \sum_{k\, =\, 0}^{n - 1} {\rm E}\big[ \left( \mathbb{h}_{w}^{k} \right)^{2} \big]\ \ \le\ \ C_{22} \left( \varsigma_{w}^{(n)} (g) \right)^{2}, \] 
for some positive constants $C_{21}$ and $C_{22}$. Hence, by Gal Koksma Theorem as in \cite{zl:14, sw:60}, we have 
\begin{eqnarray*} 
\sum_{k\, =\, 0}^{n - 1} \left( {\rm E} \left[ \left( \mathbb{h}_{w}^{k} \right)^{2} \mid \mathscr{B}_{w}^{k + 1} \right]\; -\; {\rm E} \left[ \left( \mathbb{h}_{w}^{k}\right)^{2} \right] \right) & = & O \left( \left( \varsigma_{w}^{(n)} (g) \right)\; +\; \log^{\frac{3}{2}\, +\, \epsilon} \left( \left( \varsigma_{w}^{(n)} (g) \right)^{2} \right) \right) \\ 
& = & o \left( \left( \varsigma_{w}^{(n)} (g) \right)^{1\, +\, \epsilon'} \right) \\ 
& = & o(a_{n}), 
\end{eqnarray*} 
where $\epsilon'$ is some small positive quantity, possibly less than or equal to $\epsilon$. 
\end{proof} 
\medskip 

\noindent 
\begin{proof}[Claim \eqref{cond2}] 
Here, we begin with an easy observation that $\left( \varsigma_{w}^{(n)} (g) \right)^{2} = O(n)$. Thus given $\delta > 0$, there exists a threshold $M_{\delta} \in \mathbb{N}$ and a positive real number $C_{23} > 0$ such that 
\[ \left| \frac{\varsigma_{w}^{(n)} (g)}{\sqrt{n}}\; -\; C_{23} \right|\ \ \le\ \ \delta,\ \ \ \forall n \ge M_{\delta}. \] 
Suppose $m < M_{\delta}$. Then, by the Archimedean property of the reals, we have 
\[ \varsigma_{w}^{(m)} (g)\ \ \ge\ \ \sqrt{m} D_{m}. \] 
Choosing $C_{24} = \min \big\{ C_{23} - \delta,\, D_{1},\, D_{2},\, \cdots,\, D_{N_{\delta}} \big\}$, we have 
\[ \varsigma_{w}^{(n)} (g)\ \ \ge\ \ \sqrt{n} C_{24},\ \ \ \ \forall n \in \mathbb{N}. \] 
\medskip 

\noindent 
For $r = 2$ in condition (2) in the enumerated statement of theorem \eqref{CM 2.3}, we have 
\begin{eqnarray*} 
\sum_{n\, \ge\, 0} a_{n}^{- 2} {\rm E} \left[ \left\vert \mathbb{h}_{w}^{n} \right\vert^{4} \right] & = & \sum_{n\, \ge\, 0} \left( \varsigma_{w}^{(n)} (g) \right)^{- (2\, +\, \epsilon)} {\rm E} \left[ \left\vert \mathbb{h}_{w}^{n} \right\vert^{4} \right] \\ 
& \le & C_{25}\; +\; \sum_{n\, \ge\, 1} \frac{1}{C_{24}^{2\, +\, \epsilon} n^{1\, +\, 2 \epsilon}} {\rm E} \left[ \left\vert \mathbb{h}_{w}^{n} \right\vert^{4} \right] \\ 
& \le & C_{25}\; +\; C_{26} \sum_{n\, \ge\, 1} \frac{1}{n^{1\, +\, 2 \epsilon}} \\ 
& < & \infty, 
\end{eqnarray*} 
since $\sup\limits_{n\, \ge\, 0} {\rm E} \left[ \left\vert \mathbb{h}_{w}^{n} \right\vert^{4} \right] < \infty$.  
\end{proof}
\medskip 

\noindent 
\noindent 
Finally, when we replace the standard Brownian motion with a Brownian motion $\big\{ \widetilde{\mathfrak{B}^{*}} (t) \big\}_{t\, \ge\, 0}$ such that $\widetilde{\mathfrak{B}^{*}} (t)$ has variance $t \left( \varsigma_{w} (g) \right)^{2}$, we get 
\[ Y_{w}^{n}\; -\; \widetilde{\mathfrak{B}^{*}} \left( \left( n \right) \right)\ \ =\ \ O \left( n^{\frac{1}{2}\, -\, \gamma} \right),\ \ \ \ \nu\text{-a.s., for some}\ \gamma > 0. \] 
\bigskip 

\section{Proofs of other statistical properties} 
\label{seccor}

\noindent 
In this section, we write the proofs of the other statistical properties such as the central limit theorem, weak invariance principles and the law of iterated logarithms, as mentioned in theorems \eqref{osp1} and \eqref{osp2}. 
\medskip 

\noindent 
\begin{proof}[Theorem \eqref{osp1}]

\noindent 
\begin{enumerate} 
\item {\it Proof of the central limit theorem}: Recall from the proof of theorem \eqref{asip1} that 
\begin{enumerate} 
\item[(a)] $\mathfrak{Y}_{n}\ \ \stackrel{\text{d}}{=}\ \ \Phi^{n}\ \ \forall n \ge 1$ and 
\item[(b)] $\mathfrak{Y}_{n}\ \ =\ \ \mathfrak{B} (n)\; +\; O \left( n^{\frac{1}{4}\, +\, \delta} \right),\ \ \nu$-a.e.
\end{enumerate}
Further, owing to condition (b) above, we have that for some $\epsilon > 0$, 
\[ \frac{1}{\sqrt{n}} \mathfrak{Y}_{n}\ \ =\ \ \frac{1}{\sqrt{n}} \mathfrak{B} (n)\; +\; O\left( n^{- \epsilon} \right),\ \ \nu\text{-a.e.} \] 
and therefore, $\frac{1}{\sqrt{n}} \big( \mathfrak{Y}_{n} - \mathfrak{B} (n) \big) \stackrel{{\rm p}}{\longrightarrow} 0$, {\it i.e.}, converges in probability to $0$ as $n \rightarrow \infty$. But $\frac{1}{\sqrt{n}} \mathfrak{B} (n)$ is a normal distribution with mean zero and variance $\varsigma(G)^{2}$ for all $n \geq 1$. Further, owing to condition (a), we have that 
\[ \frac{1}{\sqrt{n}} \mathfrak{Y}_{n}\ \ \stackrel{\text{d}}{=}\ \ \frac{1}{\sqrt{n}} \Phi^{n}\ \ \forall n \ge 1. \] 
Making use of both the conditions, we have 
\[ \frac{1}{\sqrt{n}} \Phi^{n}\ \ \stackrel{\text{d}}{\longrightarrow}\ \ \mathcal{N} \left( 0,\, \varsigma(G)^{2} \right), \] 
where $\mathcal{N} \left( 0,\, \varsigma^{2} \right)$ denotes the normal distribution with mean $0$ and variance $\varsigma^{2}$. The result follows since $G^{n}((w, x)) = \Phi^{n} ((w, x)) + O(1)$. 
\medskip

\item {\it Proof of the law of iterated logarithms}: If $\Phi$ in theorem \eqref{asip1} satisfies the law of iterated logarithms, then so does $G$, since 
\[ \limsup_{n\, \rightarrow\, \infty} \left[ \frac{G^{n} ((w, x)) - \Phi^{n} ((w, x))}{\varsigma(G) \sqrt{2n \log \log n}} \right]\ \ \stackrel{\text{p}}{\longrightarrow}\ \ 0\ \ \ \ \text{as}\ n \rightarrow \infty. \] 
The following lemma is the key to proving the law of iterated logarithms for the given function $\Phi$. 
\medskip 

\noindent 
\begin{lemma}
Any Brownian motion with variance $\varsigma^{2}$ satisfies the law of iterated logarithms {\it i.e.}, 
\[ \limsup_{t\, \rightarrow\, \infty} \frac{\mathfrak{B} (t) (\omega)}{\varsigma \sqrt{2t \log \log t}}\ \ =\ \ 1,\ \ \ \nu\text{-a.e.} \] 
\end{lemma}
\medskip 

\noindent 
Since we have a Brownian motion which by condition (b) in the proof of the central limit theorem is equal to $\mathfrak{Y}_{n} + O \left( n^{\frac{1}{4}\, +\, \delta} \right)$, we have by theorem \eqref{asip1} that 
\[ \limsup_{n\, \rightarrow\, \infty} \frac{\mathfrak{Y}_{n} (\omega)}{\varsigma \sqrt{2n \log \log n}}\ \ =\ \ 1,\ \ \ \nu\text{-a.e.} \] 
Since $\mathfrak{Y}_{n}$ and $\Phi^{n}$ have the same distribution with the latter having variance $\big( \varsigma (G) \big)^{2}$, we conclude that 
\[ \limsup_{n\, \rightarrow\, \infty} \frac{\Phi^{n} ((w, x))}{\varsigma(G) \sqrt{2n \log \log n}}\ \ =\ \ 1,\ \ \ \mu_{F}\text{-a.e.} \] 
\end{enumerate} 
\end{proof}
\medskip

\noindent 
We now prove theorem \eqref{osp2} to conclude the proofs of all the theorems in this paper. 
\medskip 

\noindent 
\begin{proof}[Theorem \eqref{osp2}]
The proofs of both the statements run {\it mutatis mutandis} as the proofs of their analogous statements in theorem \eqref{osp1}. Hence, we merely highlight the following for readers' convenience. 
\medskip 

\noindent
\begin{enumerate} 
\item[(a)] $Y_{w}^{n}\ \ \stackrel{\text{d}}{=}\ \ g_{w}^{n}\ \ \ \forall n \geq 1$. 
\item[(b)] $Y_{w}^{n}\ \ =\ \ \widetilde{\mathfrak{B}}^{*} (n) + O \left( n^{\frac{1}{2}\, -\, \gamma} \right),\ \ \nu$-a.e. 
\end{enumerate} 
Further, 
\[ \limsup_{n\, \rightarrow\, \infty} \frac{g_{w}^{n}(x)}{\varsigma_{w}(g) \sqrt{2n \log \log n}}\ \ =\ \ 1,\ \ \ \lambda\text{-a.e.} \] 
\end{proof}
\bigskip 

\section{Concluding remarks} 
\label{concl} 

\noindent 
As stated in the introductory section, the theorems proved in this paper are easily transferable to several other analogous settings of dynamical systems. We conclude this paper by merely pointing to some of those. 
\medskip 

\noindent 
\begin{itemize} 

\item Instead of working with the specified interval maps $T_{d} : I \longrightarrow I$ of degree $(d + 1)$ for $1 \le d \le N$ given by $T_{d} (x) = (d + 1) x \pmod 1$, one may well consider the action of any $N$ (piecewise) linear interval maps, $S_{d} : I \longrightarrow I;\ 1 \le d \le N$, each with integer degree at least $2$. The measure, in this case still remains Lebesgue. 

\item One might as well consider the simple monomial maps $P_{d};\ 1 \le d \le N$, defined on the Riemann sphere $\overline{\mathbb{C}} = \mathbb{C} \cup \{ \infty \}$ and given by $P_{d} (z) = z^{k_{d}}$ where $k_{d} \ge 2$ for all $1 \le d \le N$. We know from, say \cite{afb:91}, that the Julia set $\mathbb{J} (P_{d})$ of the monomial map $P_{d}$ is the unit circle $\mathbb{S}^{1}$ in $\mathbb{C}$. In such a case, the Julia set $\mathbb{J} (P)$ of the skew-product map $P$ appropriately defined analogous to equation \eqref{spm} is also the unit circle, as one may find from \cite{hs:00}. Owing to the Julia set $\mathbb{J} (P)$ being completely $P$-invariant, one may undertake an analogous study, as done in this paper, to the dynamics generated by the monomial maps $P_{d};\ 1 \le d \le N$ restricted on the Julia set $\mathbb{J} (P) = \mathbb{S}^{1} \subset \mathbb{C}$ and obtain analogous results, employing the Haar measure on $\mathbb{S}^{1}$. 

\item Let $R_{d};\ 1 \le d \le N$ be a collection of rational maps acting on the Riemann sphere $\overline{\mathbb{C}}$; each with degree $k_{d} \ge 2$. We suppose that the rational maps are so chosen that the Julia set $\mathbb{J} (R_{d})$ of the map $R_{d}$ is topologically connected. Then, defining the skew-product map $R$ appropriately and restricting its action on the $R$-invariant Julia set $\mathbb{J} (R)$, as one may obtain from \cite{hs:00}, it is possible to investigate analogous results for Boyd's measure, as defined in \cite{db:99}. It must be borne in mind that the Boyd's measure is a generalisation of the Lyubich's measure, as in \cite{myl:86}. 

\item Finally, we consider the dynamical system obtained by iterating certain relations on $\mathbb{C}$. The relation can be explained as the zero set of a polynomial, say $Q \in \mathbb{C}[\zeta,\, \omega]$ of a certain form such that: 
\begin{itemize} 
\item $Q \big( \cdot,\, \omega \big)$ and $Q \big( \zeta,\, \cdot \big)$ are generically multiple-valued; 
\item if $\mathscr{G}_{Q}$ denotes the biprojective completion of $\big\{ Q\, =\, 0 \big\}$ in $\overline{\mathbb{C}}\, \times\, \overline{\mathbb{C}}$, then no irreducible component of $\mathscr{G}_{Q}$ is of the form $\big\{ a \big\}\, \times\, \overline{\mathbb{C}}$ or $\overline{\mathbb{C}}\, \times\, \big\{ a \big\}$, where $a \in \overline{\mathbb{C}}$.  
\end{itemize} 
Such dynamical systems have been studied by Dinh and Sibony in \cite{ds:06} and Bharali and Sridharan in \cite{bs:16}. Upon satisfying certain technical conditions, one may study analogous dynamical and statistical properties with respect to the Dinh-Sibony measure, when the action of the holomorphic correspondence $Q \big( \zeta,\, \omega \big)$ is restricted on the support of the Dinh-Sibony measure, as defined in \cite{ds:06, bs:16}. 
\end{itemize}

\end{document}